\documentclass[reqno,12pt]{amsart}
\usepackage{amsfonts}
\usepackage{amssymb}
\usepackage{amsthm}
\usepackage{upref}
\usepackage{enumerate}

\usepackage{amscd}

\usepackage [latin1]{inputenc}

\makeatletter
\def\LaTeX{\leavevmode L\raise.42ex
    \hbox{\kern-.3em\size{\sf@size}{0pt}\selectfont A}\kern-.15em\TeX}
 \makeatother
 
\DeclareMathOperator{\clos}{clos}

\numberwithin{equation}{section}

\newtheorem{lemma}{Lemma}[section]
\newtheorem{theorem}[lemma]{Theorem} 
\newtheorem{corollary}[lemma]{Corollary}
\newtheorem{proposition}[lemma]{Proposition}

\theoremstyle{definition}

\newtheorem{remark}[lemma]{Remark}

  \newcommand{\e}{\eqref}

\newcommand{\q}{\quad}

\newcommand{\ov}{\overline}
\newcommand{\wt}{\widetilde}

\newcommand{\ti}{\tilde}

\renewcommand{\d}{\delta}

\newcommand{\la}{\langle}
\newcommand{\ra}{\rangle}

\renewcommand\Im{\operatorname{Im}}
\renewcommand\Re{\operatorname{Re}}

\newenvironment{pf}{\begin{proof}}{\end{proof}}

\def\qqq{\mathrel{\subset\mkern-15mu\lower.38ex\hbox{${\scriptscriptstyle\rightarrow}$}}}

\let\cal\mathcal

\let\Bbb\mathbb
 
\begin{document}
\title 
[ Spectral theory of Jacobi operators with increasing coefficients]
{Spectral theory of Jacobi operators with increasing coefficients. The critical case} 

\author{ D. R. Yafaev  }
\address{CNRS, IRMAR-UMR 6625, F-35000
    Rennes, France and SPGU, Univ. Nab. 7/9, Saint Petersburg, 199034 Russia}
\email{yafaev@univ-rennes1.fr}
 

\thanks {Supported by  project   Russian Science Foundation   22-11-00070}

\begin{abstract}
Spectral properties of Jacobi operators $J$  are intimately related to an asymptotic behavior of the corresponding orthogonal polynomials $P_{n}(z)$ as $n\to\infty$. We study the case where the off-diagonal  coefficients $a_{n}$   and, eventually,  diagonal coefficients $ b_{n}$   of $J$ tend to infinity  in such a way that  the ratio $\gamma_{n}:=2^{-1}b_{n} (a_{n}a_{n-1})^{-1/2} $ has a finite limit $ \gamma $.
In the case $|\gamma  |  < 1$
asymptotic formulas for
$P_{n}(z)$  generalize those for the Hermite polynomials and the corresponding Jacobi operators $J$ have  absolutely continuous  spectra covering the whole real line.  If $|\gamma  |  > 1$, then   spectra of  the operators $J$ are discrete. Our goal is to investigate  the critical case $| \gamma  |=1$  that occurs, for example,  for the Laguerre polynomials.
The formulas obtained depend crucially on the rate of growth of the coefficients $a_{n}$ (or $b_{n}$) and are qualitatively different in the cases where $a_{n}\to \infty$ faster or slower then $n$.   For the fast growth of $a_{n}$, we also have to distinguish the cases  
 $|\gamma_{n}|
 \to  1-0$  and $|\gamma_{n}|
 \to  1+0$.   Spectral properties of the corresponding Jacobi operators are  quite different in all these cases. Our approach works for an arbitrary power growth of the Jacobi coefficients.
      \end{abstract}



 \maketitle

\section{Introduction. Basic definitions }



\subsection{Jacobi operators }

We consider  Jacobi operators defined by three-diagonal matrices
\[
{\cal J} = 
\begin{pmatrix}
 b_{0}&a_{0}& 0&0&0&\cdots \\
 a_{0}&b_{1}&a_{1}&0&0&\cdots \\
  0&a_{1}&b_{2}&a_{2}&0&\cdots \\
  0&0&a_{2}&b_{3}&a_{3}&\cdots \\
  \vdots&\vdots&\vdots&\ddots&\ddots&\ddots
\end{pmatrix} 
\]
 in the canonical basis  of the space $  \ell^2 ({\Bbb Z}_{+})$.  Thus, if $u=(u_{0}, u_{1},\ldots )^\top=: (u_{n})$ is a column, then
 \[
( {\cal J} u) _{0} =  b_{0} u_{0}+ a_{0} u_{1} \q  \mbox{and} \q
( {\cal J} u) _{n} = a_{n-1} u_{n-1}+b_{n} u_{n}+ a_{n} u_{n+1} \q \mbox{for}\q n\geq 1  .
\]
 It is always supposed that $a_n >0$, $b_n=\bar{b}_n $ so that the matrix $ {\cal J} $ is symmetric and commutes with the complex conjugation. The minimal Jacobi operator $J_{\rm min}$ is  defined   by the equality $J_{\rm min} u= {\cal J} u$ on the set $\cal D \subset \ell^2 ({\Bbb Z}_{+})$ of vectors $u=  (u_{n}) $ with only a  finite number of non-zero components  $u_{n}$. The operator $J_{\rm min} $ is
symmetric  in the space $\ell^2 ({\Bbb Z}_{+})$ and $J_{\rm min}:  {\cal D}\to  {\cal D}$.  Its adjoint $J_{\rm min}^* $ coincides with the maximal  operator $J_{\rm max}$  given by the  same formula $J_{\rm max} u= {\cal J} u$ on the set $ {\cal D} (J_{\rm max})$  of all vectors $u \in \ell^2 ({\Bbb Z}_{+})$ such that ${\cal J}u \in \ell^2 ({\Bbb Z}_{+})$. 
 
 The operator $J_{\rm min} $ is bounded if and only if both sequences $a_{n}$ and $b_{n}$  are in 
$\ell^\infty ({\Bbb Z}_{+})$.
 In general, $J_{\rm min} $ may have deficiency indices $(0,0)$ (that is,  it is essentially self-adjoint) or $(1,1)$.    Its essential self-adjointness depends on a behavior of solutions to the difference equation
       \begin{equation}
 a_{n-1} F_{n-1} (z) +b_{n} F_{n } (z) + a_{n} F_{n+1} (z)= z F_n (z), \q n \geq 1. 
\label{eq:Jy}\end{equation}
Recall that the Weyl theory developed by him for differential equations can be naturally adapted to equations \e{eq:Jy} (see, e.g., \S 3 of Chapter~1 in the book \cite{AKH}  and references therein). For $\Im z\neq 0$,  equation \e{eq:Jy} always has a non-trivial solution $F_{n}(z)\in \ell^2 ({\Bbb Z}_{+})$. This solution is either unique (up to a constant factor) or all solutions of  equation \e{eq:Jy}  belong to $\ell^2 ({\Bbb Z}_{+})$.  The first instance is known as the limit point case and the second one -- as the limit circle case. It turns out that the operator $J_{\rm min} $ is essentially self-adjoint if and only if the limit point case occurs; 
then the closure $\clos J_{\rm min}$ of $J_{\rm min}$ equals $J_{\rm max}$.
In  the limit circle case, the operator $J_{\rm min} $  has deficiency indices  $(1,1)$.

It is well known that  the limit point case occurs if   $a_{n}\to\infty$ as $n\to\infty$ but  not too rapidly.  For example, the condition
   \begin{equation}
\sum_{n=0}^\infty a_{n}^{-1}=\infty
\label{eq:Carl}\end{equation}
(introduced by T.~Carleman in his book \cite{Carleman}) is sufficient for the essential self-adjointness of the operator $J_{\rm min} $. Under this condition no assumptions on the diagonal elements  $b_{n}$ are required. In general, the essential self-adjointness of  $J_{\rm min} $ is determined by a competition between sequences $a_{n}$ and $b_{n}$. For example, if $b_{n}$ are much  larger than $a_{n}$, then  $J_{\rm min} $ is close to a diagonal operator so that it is essentially self-adjoint independently of the growth of $a_{n}$.

     \subsection{Orthogonal polynomials}

 Orthogonal polynomials  $P_{n } (z)$ can be formally defined  as ``eigenvectors" of the Jacobi operators. This means that  a column  
 $$P(z)=(P_{0}(z), P_1(z),\ldots)^{\top}$$
  satisfies the equation ${\cal J} P( z )= z P(z)$ with $ z\in{\Bbb C}$ being an ``eigenvalue". This equation is equivalent to the recurrence relation 
  \begin{equation}
 a_{n-1} P_{n-1} (z) +b_{n} P_{n } (z) + a_{n} P_{n+1} (z)= z P_n (z), \q n\in{\Bbb Z}_{+}=\{0,1, 2, \ldots\}, 
\label{eq:RR}\end{equation}
 complemented by boundary conditions $P_{-1 } (z) =0$, $P_0 (z) =1$. Determining $P_{n } (z)$, $n=1,2, \ldots$, successively from \e{eq:RR}, we see that $P_{n } (z)$ is a polynomial 
with real coefficients of degree $n$: $P_{n } (z)= p_{n}z^n+\cdots$ where $p_{n}= (a_{0}a_{1}\cdots a_{n-1})^{-1}$.

   The spectra of all  self-adjoint extensions $J$ of the minimal operator $J_{\rm min} $  are simple with $e_{0} = (1,0,0,\ldots)^{\top}$ being a generating vector. Therefore  it is natural to define   the   spectral measure of $J$ by the relation $d\Xi_{J}(\lambda)=d\la E_{J}(\lambda)e_{0}, e_{0}\ra $ where  $E_{J}(\lambda)$      is the spectral family of the operator $J$ and $\la\cdot, \cdot\ra$ is the scalar product in the space  $  \ell^2 ({\Bbb Z}_{+})$.
  For all   extensions $J$ of the operator $J_{\rm min} $, the  polynomials  $P_{n}(\lambda)$ are orthogonal and normalized  in the spaces $L^2 ({\Bbb R};d\Xi_{J})$: 
      \[
\int_{-\infty}^\infty P_{n}(\lambda) P_{m}(\lambda) d\Xi_{J}(\lambda) =\d_{n,m};
\]
as usual, $\d_{n,n}=1$ and $\d_{n,m}=0$ for $n\neq m$. We always consider normalized  polynomials  $P_{n}(\lambda)$. They  are often called orthonormal. If the operator $J_{\rm min} $ is essentially self-adjoint and $J=\clos J_{\rm min} $, we write $d\Xi(\lambda)$ instead of  $d\Xi_{J}(\lambda)$.
 
 It is useful to keep in mind the following   elementary observation.
 
   \begin{proposition}\label{refl}
       If  a sequence  $F_{n}(z)$   satisfies equation \e{eq:Jy}, then 
\[
F^\sharp _{n}(z)= (-1)^{n} F_{n}(-z)
\]
satisfies the same equation with
the Jacobi coefficients $( a^\sharp_{n},  b^\sharp_{n}) = (a_{n}, - b_{n})$. In particular, $P^\sharp_{n}(z)= (-1)^{n} P_{n}(-z)$ are the orthonormal polynomials
for the coefficients   $( a^\sharp_{n},  b^\sharp_{n})$.  In the limit point case, 
if $  J^\sharp $ is the Jacobi operator in the space $  \ell^2 ({\Bbb Z}_{+})$ with matrix elements $( a^\sharp_{n},  b^\sharp_{n})$, then $J^\sharp =-{\cal U}^* J {\cal U}$ where the unitary operator ${\cal U}$ is defined by $({\cal U} F)_{n}=(-1)^n F_{n}$ for $n\in{\Bbb Z}_{+}$.   The corresponding spectral measures are linked by the relation
$d\Xi^\sharp(\lambda)= d \Xi (-\lambda)$. In particular, if $b_{n}=0$ for all $n$, then the operators $J$ and $-J$ are unitarily equivalent.  
\end{proposition}
 
 The comprehensive presentation of the results  described shortly above  can be found in  the books \cite{AKH, Chihara,   Sz} and the surveys \cite{Lub,Simon, Tot,  Asshe}.
     
   \subsection{Asymptotic results}

We study   the case $a_{n}\to\infty$  as $n\to\infty$ and are interested in the asymptotic behavior of the polynomials  $P_{n } (z)$ as $n\to\infty$.  The  condition  $a_{n}\to\infty$ is fulfilled for  the Hermite polynomials where the Jacobi coefficients are
 \begin{equation}
 a_{n}=\sqrt{(n+1)/2}  \q\mbox{and}\q    b_{n} =0
 \label{eq:Her}\end{equation}
and the Laguerre polynomials  
$L_{n}^{(p)}(z) $ where
     \begin{equation}
    a_{n} = \sqrt{(n+1)(n+1+p)} \q\mbox{and}\q     b_{n} =  2n+p+1, \q p>-1. 
\label{eq:Lag}\end{equation}
In the general case there are two essentially different approaches to this problem.  The first one derives asymptotic formulas for  $P_{n } (z)$ from the spectral measure $d\Xi(\lambda)$, and the second proceeds directly from the coefficients $a_{n}$, $b_{n}$.  The first method goes back to  S.~Bernstein 
(see his pioneering paper \cite{Bern} or  Theorem~12.1.4 in the G.~Szeg\H{o} book~\cite{Sz})  who obtained formulas generalizing those for the Jacobi polynomials. In terms of the coefficients $a_{n}$, $b_{n}$, the assumptions of  \cite{Bern}  correspond to the conditions
     \begin{equation}
a_{n} \to a_{\infty}>0, \q b_{n} \to  0\q \mbox{as}\q n\to\infty.
\label{eq:Nev}\end{equation}
   Generalizations of the asymptotic formulas for the Hermite polynomials are known as the Plancherel-Rotach formulas.


A study of an asymptotic behavior of the orthonormal  polynomials for
given  coefficients $a_{n}$, $b_{n}$  was initiated by P.~Nevai in his book \cite{Nev}. He (see also the papers \cite{Mate} and  \cite{Va-As}) investigated  the case of stabilizing coefficients  satisfying condition \e{eq:Nev}, but, in contrast to  \cite{Bern},  the results of  \cite{Nev, Mate, Va-As} were stated directly in terms of the Jacobi coefficients.  
The case of the coefficients  $a_{n}\to \infty$     was later  studied in \cite{Jan-Nab} by J.~Janas and  S.~Naboko and in \cite{Apt} by A.~Aptekarev and J.~Geronimo. It was assumed in these papers  that  there exists a finite   limit
 \begin{equation}
   \frac{b_{n}}{2\sqrt{a_{n-1}a_{n}} }=: \gamma_{n}\to  \gamma ,  \q n\to\infty,
\label{eq:Gr}\end{equation}
  where $| \gamma | < 1$ so that $b_{n}$  are relatively small compared to $a_{n}$. The   Carleman condition
\e{eq:Carl}  was also required.
 The famous example of this type is given by the Hermite coefficients \e{eq:Her}.
  In the general case the results are qualitatively similar to this particular case. Asymptotics of $P_{n} (\lambda)$ are oscillating for $\lambda\in \Bbb R$ and $P_{n} (z)$ exponentially grow as $n\to\infty$ if $\Im z\neq 0$.  Spectra of the operators $J$ are absolutely continuous and fill the whole real axis. 
    If  \e{eq:Gr}  is satisfied with $| \gamma | >1$, then    diagonal elements $b_{n}$ dominate off-diagonal elements $a_{n}$. This ensures that the spectra of such operators $J$ are  discrete.  Note (see, e.g., \cite{inf}) that algebraic structures of
     asymptotic formulas for the orthonormal  polynomials    are quite similar in  the cases $| \gamma | <1$ and $| \gamma | >1$, but in the second case    $P_{n} (z)$ exponentially grow as $n\to\infty$ even for $z\in\Bbb R$ (unless $z$ is an eigenvalue of $J$).

 The case of rapidly increasing coefficients $a_{n}$ when the Carleman condition \e{eq:Carl} is violated,  so that
  \[
  \sum_{n =0}^{\infty} a_{n}^{-1} <\infty ,
 \]
 was investigated in a recent paper \cite{nCarl} where it was also assumed   that $|\gamma | \neq 1$. 
  Astonishingly, the asymptotics of the orthogonal polynomials in this a priori highly singular case is particularly simple and general.

 \subsection{Critical case}

In  the critical case   $|\gamma| = 1$, the coefficients $a_{n}$ and $b_{n}$ are of the same order and asymptotic formulas for $P_{n}(z)$ are determined by details of their behavior as $n\to\infty$.


Thus, one has to require   assumptions on the coefficients $a_{n}$ and $b_{n}$ more specific  compared to \e{eq:Gr}.  
To make our presentation as simple as possible, we assume that, asymptotically, 
  \begin{equation}
  a_{n}=  n^\sigma (1+  \alpha n^{-1}+ O (n^{-2}) ),\q n\to\infty, 
\label{eq:ASa}\end{equation}
and
  \begin{equation}
  b_{n}=  2 \gamma n^\sigma  (1+  \beta n^{-1}+ O (n^{-2}) ),\q n\to\infty,  
\label{eq:ASb}\end{equation}
for some   $\alpha,\beta, \gamma\in{\Bbb R}$ and\footnote{The case $\sigma>3/2$ was considered earlier in \cite{crit}} $\sigma>0$. Thus, the  operators with periodically modulated coefficients (see, e.g., \cite{DN} and references therein) are out of the scope of this paper. The critical case is distinguished  by the condition $|\gamma|=1$. In view of Proposition~\ref{refl}  the results for $\gamma=1$ and $\gamma=-1$ are equivalent. It turns out that the asymptotic formulas for $P_{n}(z)$  depend crucially on the parameter
 \begin{equation}
\tau= 2 \beta -2\alpha + \sigma .  
\label{eq:BX3}\end{equation}
Roughly speaking, the cases $\tau<0$ (or $\tau>0$)  correspond to dominating off-diagonal $a_{n}$ (resp., diagonal $b_{n}$) Jacobi coefficients.

  All the results of this paper   can be extended to a more general situation where the terms $\alpha n^{-1}$ and $\beta n^{-1}$ in
  \e{eq:ASa},  \e{eq:ASb} are replaced by $\alpha n^{-p}$ and $\beta n^{-p}$ for some $p\in (0,2)$ and the error term $O(n^{-2})$ is replaced by
  $O(n^{-r})$  for $r>\max\{1,p\}$.

  The classical example where the critical case occurs is given by  the Laguerre   coefficients \e{eq:Lag}.   In this case, we have  $\gamma=1$, $\sigma=1$ and $\alpha=1+ p/2$, $\beta= (1+p)/2$ so that $\tau=0$.
    The corresponding Jacobi operators $J=J^{(p)}$ have absolutely continuous  spectra
coinciding with $[0,\infty)$.  Another  example is given by the Jacobi operators describing birth and death processes investigated in \cite{Jan-Nab2} and  \cite{Mosz}. The recurrence coefficients of such operators are rather close to \e{eq:Lag} so that spectral and asymptotic results for these two classes of operators are similar.

  Probably, a study of   Jacobi operators   in  the critical case was initiatiated by J.~Dombrowsi and S.~Pedersen in  the   papers \cite{Domb-P, Domb}  where spectral properties of such operators were investigated  under sufficiently general assumptions on the coefficients $a_{n}$ and $b_{n}$. Asymptotics of the orthogonal polynomials in this  situation was studied by  J.~Janas, S.~Naboko and E.~Sheronova in
the pioneering paper \cite{Jan-Nab-Sh}. They accepted conditions \e{eq:ASa},  \e{eq:ASb}  with $\sigma\in (1/2, 2/3)$, $\alpha=\beta=0$ and studied equation \e{eq:Jy} for real $z=\lambda$.  Both  
  oscillating for $\lambda>0$ (if $\gamma=1$) and exponentially growing   (or  decaying)  for $\lambda<0$  (if $\gamma=1$)  asymptotics of solutions of equation   \e{eq:Jy}  were investigated in \cite{Jan-Nab-Sh}.  The results of this  paper imply that  positive spectra of the operators $J$ are absolutely continuous and  negative spectra are discrete.  Recently the results of  \cite{Jan-Nab-Sh} were generalized and supplemented in  \cite{Na-Si}  by some ideas of \cite{Apt}  -- see Remark~\ref{Na-Si} below.
  
  We note also the paper \cite{Sah} by J.~Sahbani where  interesting spectral  results were obtained avoiding a study of asymptotics of the orthogonal polynomials. The paper \cite{Sah}  relies on the Mourre method.

In the non-critical case  $|\gamma|\neq 1$, asymptotic formulas are qualitatively different for  $\sigma\leq 1$  when the Carleman condition is satisfied and for $\sigma>1 $ when the Carleman condition fails. In the  critical case the borderline is
  $\sigma=3/2$.  The case of rapidly increasing coefficients where $\sigma>3/2$ was studied in \cite{crit}. For such $\sigma$, the limit circle case is realized (if $\tau<0$) and the corresponding Jacobi operators have discrete spectra.
  
  Our goal is to consistently study
  the regular critical case
  where $|\gamma|=1$ and $\sigma \leq 3/2$. Then the Jacobi operator $J_{\rm min}$ is essentially self-adjoint, even if the Carleman condition   \e{eq:Carl}  fails.
  Its spectral properties turn out to be qualitatively different in the  cases $\sigma\in (0,1)$, $\sigma=1$ and $\sigma\in (1,3/2]$.  Moreover, for $\sigma\in (1,3/2]$ the answers depend crucially on the sign of the parameter  $\tau$ defined by  \e{eq:BX3}. 
In all cases, our asymptotic formulas are constructed in terms of the sequence
  \begin{equation}
t_{n} (z)=-\tau n^{-1}+  z n^{-\sigma}  .
\label{eq:ME1}\end{equation}

  Note that the critical situation studied here is morally similar to a threshold behavior of orthogonal polynomials for  case  \e{eq:Nev}. For such coefficients,  the role of \e{eq:Gr} is played (see \cite{Nev,Mate,JLR}) by the relation 
   \[
 \lim_{n\to\infty}  \frac{b_{n}-\lambda}{2 a_{n}} =-\frac{ \lambda}{2 a_\infty} .
\]
Since the essential spectrum of the operator $J$ is now $[-2 a_\infty,2 a_\infty]$, the values $\lambda=\pm 2 a_\infty$ are the threshold values of the spectral parameter $\lambda$.  The parameter $-\lambda/ (2 a_\infty)$ plays the role of $\gamma $ so that the cases $|\gamma |<1$ (resp., $|\gamma |>1$)  correspond to $\lambda$ lying inside the essential spectrum of   $J$  (resp., outside of it).

    \subsection{Scheme of the approach}
   
    
    We use the traditional  approach developed for
    differential equations
 \begin{equation}
 - (a(x) f ' (x, z) )'+ b(x) f  (x, z)= z f  (x, z), \q  x>0,  \q a(x) > 0.
\label{eq:Schr}\end{equation}
  To a large extent, $x$, $a(x)$ and $b(x) $ in \e{eq:Schr} play the roles of the parameters $n$, $a_n$ and $b_{n}$ in the Jacobi equation \e{eq:Jy}. The regular solution $\psi(x,z)$ of the differential equation \e{eq:Schr}  is distinguished by the conditions
 \[
\psi (0, z) =0,\q   \psi' (0, z) =1.
\]
It plays the role of the polynomial solution $P_n (z)$ of equation \e{eq:Jy} fixed by the conditions
$P_{-1} (z)=0$, $P_0 (z)=1$.

A study of an asymptotics of the regular solution $\psi(x,z)$  relies on a construction of   special solutions of the differential equation  \e{eq:Schr}   distinguished by their asymptotics as $x\to\infty$. For example, in the case $a(x)=1$, $b\in L^1 ({\Bbb R}_{+})$,  equation  \e{eq:Schr} has a solution $f  (x, z)$,  known as the Jost solution, behaving like $e^{i\sqrt{z}x}$, $  \Im \sqrt{z} \geq 0$, as $x\to\infty$.
Under fairly general assumptions  
equation \e{eq:Schr} has a solution $f(x,z)$ (we also  call it the Jost solution) whose asymptotics is given by  the classical Liouville-Green formula (see   Chapter~6 of the book \cite{Olver})
 \begin{equation}
 f  (x, z)\sim   {\cal G} (x,z)^{-1/2}  \exp  \Big(i \int_{x_{0}} ^x  {\cal G} (y,z) dy\Big)=: {\cal A}  (x, z)
\label{eq:Ans}\end{equation}
as $x\to\infty$.  Here $x_{0}$ is some fixed number and
\[
{\cal G} (x,z)= \sqrt{\frac{z-b(x)}{a(x)} } , \q\Im {\cal G} (x,z)  \geq 0.
\]
 Note that the function $ {\cal A}  (x, z)$ (the Ansatz for the  Jost solution $f  (x, z)$) satisfies equation \e{eq:Schr} with a sufficiently good accuracy.   
 
   For real $\lambda$ in the absolutely continuous spectrum of the operator 
   \[
   -\frac{d}{dx} \big(a(x)\frac{d}{dx} \big)+b(x), 
   \]
   the regular solution $\psi (x,\lambda)$  of  \e{eq:Schr}  is a linear combination of the Jost solutions $f(x,\lambda+i 0)$ and $f(x,\lambda-i0)$ which yields asymptotics of $\psi (x,\lambda)$  as $x\to \infty$.  For example, in the case $a(x)=1$, $b\in L^1 ({\Bbb R}_{+})$ and $\lambda>0$, one has
  \[
\psi(x,\lambda) =  \kappa(\lambda) \sin (\sqrt{\lambda}  x+ \eta (\lambda)) + o(1), \q x\to\infty,
\]
where $\kappa(\lambda)$ and $\eta(\lambda)$  are known as the scattering (or limit) amplitude and phase, respectively.
  If $\Im z\neq 0$, then one additionally constructs, by an explicit formula,  a solution $g(x,z)$ of \e{eq:Schr} exponentially growing as $x\to\infty$. This yields 	asymptotics of  $\psi(x,z)$ for $\Im z\neq 0$. 
  

 An analogy between the equations \e{eq:Jy} and \e{eq:Schr} is of course very well known. However it seems to be never consistently exploited before. In particular, the papers cited above use also specific methods of difference  equations. For example,   the absolute continuity of the spectrum is often deduced from the subordinacy theory, the asymptotics of the orthonormal polynomials are calculated by studying infinite products of transfer matrices, etc.  Some of these tools are quite ingenious, but, in the author's opinion, the standard approach of differential equations works perfectly well and allows one to study an asymptotic behavior of orthonormal polynomials in a very direct way. It permits  an arbitrary growth of the coefficients $a_{n}$ and $b_{n}$ (all values of $\sigma$ in formulas \e{eq:ASa}, \e{eq:ASb}) and naturally leads   to a variety of new results, for example, to a construction of the resolvents of Jacobi operators and to  the limiting absorption principle. For Jacobi operators with increasing coefficients, this approach was already used in the non-critical case  $|\gamma |\neq 1$  in \cite{inf}.

We are applying the same scheme to the  regular critical  case  when conditions  \e{eq:ASa}  and \e{eq:ASb}  are satisfied with $\sigma\leq 3/2$ and  $| \gamma |=1$ in \e{eq:ASb}. Under these assumptions the limit point case occurs although for $\sigma>1$ the Carleman condition \e{eq:Carl}  is violated. 

  Let us briefly describe the main steps of our approach. In the non-critical case $| \gamma |\neq1$, it was presented in \cite{inf}.

A. 
First, we  distinguish solutions (the Jost solutions) $f_{n}  (z)$ of the  difference  equation \e{eq:Jy} 
 by their asymptotics as $n\to\infty$. This requires a construction of an Ansatz ${\cal A}_{n} (z)$   for   the Jost solutions   such that
  the relative remainder
\begin{equation}
 {\bf r}_{n} (z)  : =( \sqrt{ a_{n-1} a_{n}} {\cal A}_{n} (z) )^{-1} \big(a_{n-1} {\cal A}_{n-1} (z)  + (b_{n}-z){\cal A}_{n}(z)   + a_{n} {\cal A}_{n+1} (z) \big) 
\label{eq:Grr}\end{equation}
belongs at least to the space $ \ell^1 ({\Bbb Z}_{+})$.

 
  B.
We seek ${\cal A}_{n}(z) $ in the form
    \begin{equation}
{\cal A}_{n}(  z )  =  (-\gamma)^n n^{- \rho} e^{  i \varphi_{n} (\gamma z)} , \q \gamma=\pm 1,
\label{eq:Gr1}\end{equation} 
where the power $\rho$ in the amplitude and the phases $\varphi_{n} $ are determined by the coefficients $a_{n}$, $b_{n}$. Post factum,
    ${\cal A}_{n}(z) $  turns out to be the leading term of the asymptotics  of $f_{n}(  z ) $ as $n\to\infty$:
 \begin{equation}
f_{n} (z)={\cal A}_{n}(  z )  (1+ o(1)).
\label{eq:Gr1x}\end{equation}

  Actually,  the Ans\"atzen we use  are only distantly similar to the Liouville-Green Ansatz    \e{eq:Ans}.  
  On the other hand, for  $\sigma=1$, relation  \e{eq:Gr1}  is close to formulas of the Birkhoff-Adams method significantly polished in \cite{W-L} (see also Theorem~8.36 in the book \cite{El}).

C.
Then we make   a multiplicative change of variables
 \begin{equation}
   f_{n} (z)= {\cal A}_{n} (z) u_{n} (z ) 
      \label{eq:Jost}\end{equation} 
which permits us to reduce   the Jacobi equation \e{eq:Jy} for  $  f_{n} (z)$  to a Volterra ``integral" equation for the sequence $u_{n} (z)$.   
This equation  depends of course on the parameters $a_{n}$, $b_{n}$. In particular, for $\sigma>1$, it is  qualitatively   different in the cases $\tau < 0$ and  $ \tau > 0$. However in all cases the Volterra equation for $u_{n} (z)$ is standardly solved    by iterations  which allows us to prove  that it has a solution such that
$ u_{n} (z)\to 1$ as $n\to\infty$. Then the Jost solutions $ f_{n} (z)$  are defined by formula \e{eq:Jost}.


   
   D. 
   To find an asymptotics of all solutions of   the Jacobi equation \e{eq:Jy} and, in particular, of the orthonormal polynomials $P_{n} (z)$, we have to construct a solution linearly independent with $f_{n}  (z)$.
   If a real $z=\lambda$  belongs to the absolutely continuous spectrum of the operator $J$, then the solutions $ f_{n}( \lambda + i0)$
   and its complex conjugate  $ f_{n}( \lambda - i0)$
 are  linearly independent.
For regular points $z$,   a solution $ g_{n} (z ) $ of \e{eq:Jy} linearly independent with $ f_{n} (z ) $ is constructed      (see, e.g., Theorem~2.2 in \cite{inf}) by an explicit formula
    \begin{equation} 
g_{n} (z ) = f_{n} (z )\sum_{m=n_{0}}^n (a_{m-1} f_{m-1}(z) f_{m}(z))^{-1},\q n\geq n_{0} ,
\label{eq:GEg}\end{equation}
where $n_{0}=n_{0}(z)$ is a sufficiently large number. It follows from \e{eq:Gr1}, \e{eq:Gr1x} that
this solution grows exponentially (for $\sigma<3/2$) as $n\to\infty$:
     \begin{equation}
g_{n}(  z )  =  i \varkappa(z)(-\gamma)^{n+1} n^{- \rho} e^{-i\varphi_{n} (\gamma z)}\big(1 + o( 1)\big) ;
\label{eq:A2P3}\end{equation} 
the factor $ \varkappa (z)$ here  is given by equality  \e{eq:wr}, but it is inessential in \e{eq:A2P3}.
Since $g_{n} (  z )$ is linearly independent with $f_{n} (  z )$, the polynomials $P_{n} (z) $  are linear combinations of $ f_{n} (  z )$ and $g_{n} (  z )$ 
which  yields asymptotics of $P_{n} (z) $.  

 E.
  Our results on the Jost solutions $ f_{n}  (z) $ allow us to determine the spectral structure of the operator $J$ and to construct its resolvent $R(z)$.
    At  the same time, we obtain the limiting absorption principle for the operator $J$ stating that matrix elements  of its resolvent $R(z)$, that is   the scalar products $\langle R(z ) u,v \rangle $, $\Im z \neq 0$,    are continuous functions of $z$ up to the absolutely continuous spectrum of the operator $J$ if elements $u$ and $v$ belong to a suitable dense subset of $\ell^2 ({\Bbb Z}_{+})$.


 All these steps, except possibly the construction of the exponentially growing solution $g_{n} (  z )$, are rather standard. No more specific tools are required in the   problem considered.
 
 Actually, the scheme described above works virtually in all asymptotic problems in the limit point case, both for difference and differential operators. In the limit circle case, some modifications are required; see \cite{nCarl, crit}. The important differences are that, in the limit circle case,  one has two natural Ans\"atzen ${\cal A}_{n}^{(\pm)}= n^{-\rho}e^{\pm i\varphi_{n}}$ where $\varphi_{n}=\bar{\varphi}_{n}$ does not depend on the spectral parameter $z\in{\Bbb C}$ and $\rho>1/2$ so that  ${\cal A}_{n}^{(\pm)}\in {\ell}^2({\Bbb Z}_{+})$.

To emphasize the analogy between differential and difference equations, we often  use the  ``continuous" terminology (Volterra  integral equations, integration by parts, etc.) for sequences labelled by the discrete variable $n$.

Our plan is the following. The main results of the paper are stated  in Sect.~2. In Sect.~3,  we define the  number $\rho$ and the phases $\varphi_{n}$ in formula  \e{eq:Gr1} for  the Ansatz  ${\cal A}_{n} (z)$ and check an estimate
 \begin{equation}
{\bf r}_{n} (z) =O (n^{ -\d}), \q n\to\infty,
\label{eq:rem}\end{equation}  
with an appropriate $\d=\d(\rho)>1$ for remainder  \e{eq:Grr}.
A Volterra integral equation for $u_{n} (z)$ is introduced and investigated in Sect.~4. This leads to a construction of the Jost solutions $f_{n}(z)$ in Sect.~5. In this section,  the proofs of Theorems~\ref{GSS-x}, \ref{GSS-y}  and \ref{GSS-z} are concluded.
Asymptotics of the orthonormal  polynomials $P_{n} (z) $  are found in Sect.~6. The results  for regular points  $z$  and for $z$ in the absolutely continuous spectrum of the Jacobi operator $J$  are stated in Theorems~\ref{P2}  and \ref{P1}, respectively. The results on spectral properties of the Jacobi operators are
collected in Theorem~\ref{S-Adj}. Its proof is given  in Sect.~7.  
   
 \section{Main results}
 
  Our goal is to study the critical case when assumptions  \e{eq:ASa} and \e{eq:ASb} are satisfied with $| \gamma | =1$.  
  In proofs, we may suppose that $\gamma=1$. The results for $\gamma=-1$ then follow
from Proposition~\ref{refl}.
  
     The results   stated below  crucially depend on the values of $\sigma$ and $\tau$. In the cases $\sigma\in (1, 3/2]$ ($\sigma\in (0,1)$)  the first (resp., the second) term in \e{eq:ME1} is dominating  so that  the asymptotic  formulas are qualitatively different in  these cases.

     \subsection{Jost solutions}
     
     Our approach relies on a study of   solutions $f_{n} (z)$  of the Jacobi equation \e{eq:Jy} 
     distinguished
by their behavior for $n\to\infty$.  Actually, we determine the sequences  $f_{n} (z)$  by their asymptotics
   \begin{equation}
f_{n}(  z )  =  (-\gamma)^n n^{- \rho} e^{  i \varphi_{n} (\gamma z)} \big(1 + o( 1)\big) , \q n\to \infty.
\label{eq:A22G+}\end{equation} 
Here
 \begin{equation} 
\rho=  \begin{cases} \sigma /2  - 1/4 \q \mbox{for}\q\sigma \geq  1
\\ \sigma/4 \q \mbox{for}\q\sigma \leq  1
\end{cases}
\label{eq:rho}\end{equation}
(observe that $\rho$ takes the critical value $\rho=1/2$ for the critical value $\sigma=3/2$)
and
 \begin{equation} 
\varphi_{n}(z)=\sum_{m=0}^n \theta_{m} (z).
\label{eq:sum}\end{equation}
The terms $\theta_{n} (z)$ will be defined  by explicit formulas below in this subsection. Note that
    \begin{equation}
 \Im \theta_{n} (z)\geq 0.
\label{eq:the}\end{equation}
  By an analogy with differential equations, it is natural to use the term ``Jost solutions" for $f_{n} (z)$. In the situation we consider, formula \e{eq:A22G+} plays the role of the Liouville-Green formula \e{eq:Ans}.  Observe that, for an arbitrary constant $C(z)$, the sequence $C(z) f_{n} (z)$ can be also taken for the Jost solution. In particular, a finite number of terms in  equality \e{eq:sum} is inessential.
 

 We denote $\Pi={\Bbb C}\setminus{\Bbb R}$ and $\Pi_{0}={\Bbb C}\setminus{\Bbb R}_{+}$.
  The sequence $t_{n} (z)$ is   given by formula \e{eq:ME1}  where $\tau$ is number \e{eq:BX3}. The analytic function $\sqrt{t}$ is defined on $\Pi_{0}$ and $\Im \sqrt{t}>0$ for $t\in\Pi_{0}$.    Below $C$, sometimes with indices,  and $c$ are different positive constants whose precise values are of no importance.

 
  We state the results about the Jost solutions $f_{n} (z)$  separately for the cases $\sigma\in (1,3/2]$, $\sigma\in (0, 1)$ and $\sigma = 1$.  Let us start with the case $\sigma > 1$.
  

\begin{theorem}\label{GSS-x} 
       Let  assumptions  \e{eq:ASa},  \e{eq:ASb} with $| \gamma |=1$  and $\sigma\in (1,3/2 ]$ be satisfied.  
        Set $\rho= \sigma /2  - 1/4$, 
  \begin{equation}
\theta_{n}(z) =   \sqrt{t_{n} (z)}  
\label{eq:bt1}\end{equation}
 and let $\varphi_{n}(z)$ be sum  \e{eq:sum}.

If $\tau<0$, then
    for every  $z\in \clos\Pi$    equation \e{eq:Jy} has a solution $  f_{n}( z ) $ with asymptotics \e{eq:A22G+}.
For all $n\in {\Bbb Z}_{+}$, the functions $f_{n}( z )$ are analytic in $\Pi$ and are continuous up to the cut along the real axis.  

If $\tau > 0$, then asymptotic formula \e{eq:A22G+} is true for all $z\in {\Bbb C}$. In this case  the functions $f_{n}( z )$ are analytic in the whole complex plane ${\Bbb C} $.

For all $\tau\neq 0$,   formula \e{eq:A22G+} is uniform in $z$ from compact subsets    of $  \Bbb C $.
\end{theorem}

We emphasize that the asymptotic behavior of the solutions $f_{n} (z)$ as $n\to\infty$ 
is  drastically different for small diagonal elements $b_{n}$ when $\tau<0$ and
for large $b_{n}$ when $\tau>0$ -- cf. formulas \e{eq:Af} and \e{eq:Af+}, below. This manifests itself in spectral properties of the corresponding Jacobi operators $J$ -- see part~$1^0$  of Theorem~\ref{S-Adj}.

\begin{remark}\label{GSS-x1} 
Formula \e{eq:A22G+}  is true  for all $\sigma>3/2$, but in this case it can be simplified by setting $z=0$ in the right-hand side of \e{eq:A22G+}. Thus, the leading term of the asymptotics of $f_{n}(z)$ does not depend on $z\in \Bbb C $ and the power $\rho>1/2$ so that $f_{n}(z)\in \ell^2 ({\Bbb Z}_{+})$. This leads to important spectral consequences:  for $\sigma>3/2$ the deficiency indices of the minimal Jacobi operator $J_{\rm min}$ are $(1,1)$, and the spectra of all its self-adjoint extensions are discrete. The case $\sigma>3/2$  was investigated  in \cite{crit}.
\end{remark}
 
  Let us pass to the case $\sigma<1$.  
The phases $\theta_{n} (z)$ are again defined by formula \e{eq:bt1} for $\sigma> 2/3$, but
  their construction 
   is more complicated 
   for $\sigma \leq 2/3$.  Let us set
   \begin{equation}
T_{n}(z) =   t_{n} (z)  + \sum_{l=2}^L  p_{l} t_{n}^l(z) 
\label{eq:te}\end{equation}
where    a sufficiently large $L$ depends on $\sigma$ and  the real numbers $p_{l}$ are defined in Lemma~\ref{ME}. In particular, $T_{n}(z)= t_{n} (z) $ for $ \sigma>2/3$. Given $T_{n} (z)$, the phases $\theta_{n} (z)$ are defined by the formula
    \begin{equation}
\theta_{n}(z) =    \sqrt{T_{n} (z)}  
\label{eq:bt1-}\end{equation}
playing the role of \e{eq:bt1}.   It is easy to show (see Remark~\ref{the}, for details) that  $T_{n} (z)\in \Pi_{0}$; thus, $\theta_{n}(z) $ are correctly defined.


 \begin{theorem}\label{GSS-y}
       Let  assumptions  \e{eq:ASa},  \e{eq:ASb} with $| \gamma |=1$  and $\sigma\in (0,1)$ be satisfied.  Set 
  $\rho= \sigma /4$ and define the functions $\theta_{n}(z)$ by formulas \e{eq:te}, \e{eq:bt1-}. Let   $\varphi_{n}(z)$ be sum \e{eq:sum}. Then for  every $z\neq 0$ such that $  z\in \gamma \clos \Pi_{0}$,     equation \e{eq:Jy} has a solution $  f_{n}( z ) $ with asymptotics
\e{eq:A22G+}.
For all $n\in {\Bbb Z}_{+}$, the functions $f_{n}( z )$ are analytic in $z\in \gamma\Pi_{0} $ and are continuous up to the cut along the   half-axis $\gamma{\Bbb R}_{+}$, with a possible exception of the boundary point $z=0$.  
         \end{theorem}

   In   the intermediary  case $\sigma =1$, the definition of the phases $\theta_{n} (z)$ is particularly explicit and the construction of the Jost solutions is simpler than for $\sigma \neq 1$.
   
 \begin{theorem}\label{GSS-z}
       Let  assumptions  \e{eq:ASa},  \e{eq:ASb} with $| \gamma |=1$  and $\sigma=1$ be satisfied.  
       Set 
  $\rho= 1 /4$,  define the functions $\theta_{n}(z)$ by the formula  
     \[
   \theta_{n} (z)=\sqrt{-\tau + \gamma z } n^{-1/2}  ,
\]
  and let   $\varphi_{n}(z)$ be sum \e{eq:sum}. 
     Then for  every $z$ such that $z\in \gamma (\tau + \clos \Pi_{0})$, $z\neq \gamma\tau$,   equation \e{eq:Jy} has a solution $ f_{n}( z )$ with asymptotics \e{eq:A22G+}.  
For all $n\in {\Bbb Z}_{+}$, the functions $f_{n}( z )$ are analytic in $z\in \gamma (\tau+ \Pi_{0})$  and are continuous up to the cut along the   half-axis $\gamma(\tau+ {\Bbb R}_{+})$, with a possible exception of the boundary  point $\gamma\tau$.  
              \end{theorem}
              
 We emphasize that in the case $\sigma \leq 1$ the condition $\tau\neq 0$ is not required. 
 
  It is convenient to introduce a notation
 \begin{equation} 
{\cal S}=
\begin{cases}
 {\Bbb R} \q &\mbox{if}\q\sigma \in (1, 3/2],   \, \tau <0
\\ 
\emptyset \q &\mbox{if}\q  \sigma \in (1, 3/2],   \, \tau >0
 \\ 
\gamma (0, \infty) \q &\mbox{if}\q\sigma \in (0,1) 
 \\ 
\gamma (\tau, \infty) \q &\mbox{if}\q\sigma =1 .
\end{cases}
\label{eq:SS}\end{equation}
We will see in Sect.~2.4 that the spectrum of the operator $J$ is absolutely continuous on the closed interval $\clos{\cal S}$,  and it may be only discrete on  ${\Bbb R}\setminus\clos{\cal S}$. Note that Theorems~\ref{GSS-x}, \ref{GSS-y} and \ref{GSS-z}  give asymptotic formulas for the Jost solutions $f_{n} (z)$ for all $z$ in the complex plane with the cut along ${\cal S}$, except the thresholds in the absolutely continuous spectrum ($z=0$ if $\sigma\in (0,1)$ and $z=\gamma\tau$ if $\sigma=1$).
For $\lambda\in {\cal S}$,  equation   \e{eq:Jy} has two linearly independent solutions $f_{n} (\lambda+ i0)$ and its complex conjugate 
         \[
f_{n}( \lambda-i0)=\ov{f_n( \lambda+i0)}.
 \]

 Under the assumptions of any of these theorems the solution $f_{n} (z)$ of  equation     \e{eq:Jy}  is determined  essentially uniquely  by its asymptotics 
  \e{eq:A22G+}. This is discussed in Sect.~5.5 (see Propositions~\ref{uniq}, \ref{uniq+} and Remark~\ref{uniq-r}).
  
  Note that the values of $u_{m-1}$ and $u_{m }$ for some $m\in{\Bbb Z}_{+}$ determine the whole sequence $u_{n}$ satisfying the difference equation \e{eq:Jy}.   Therefore  it suffices to construct sequences $f_{n} (z)$ for sufficiently large $n$ only. Then they are extended  to all $n$ as solutions of equation \e{eq:Jy}.     
  
 We also mention that $f^\sharp _{n}(z)= (-1)^{n} f_{n}(-z)$ is the Jost solution for the Jacobi equation \e{eq:Jy} with
the  coefficients $( a^\sharp_{n},  b^\sharp_{n}) = (a_{n}, - b_{n})$.

  \subsection{Asymptotics at infinity}
  
  Here we find explicit asymptotic formulas for the phases $\theta_{n} (z)$ and then for their sums $\varphi_{n}(z)$ as $n\to\infty$. These formulas  depend crucially on the 
  values of the parameters $\sigma$ and $ \tau$.
  
  Suppose first that $\sigma\in (1, 3/2]$  and that $z \in \clos\Pi$ for $\tau<0$ and $z\in \Bbb C$ for $\tau> 0$. Then 
   the term $-\tau  n^{-1} $ is dominating in \e{eq:ME1}  so that according to  definition  \e{eq:bt1}
      \begin{equation}
 \theta_{n} (z)=  n^{-1/2} \sqrt{ |\tau |+  z n^{1-\sigma}}=\pm  \sqrt{ |\tau| } n^{-1/2} \pm \frac{ z } { 2  \sqrt{ |\tau| } }n^{1/2-\sigma} + O(n^{3/2-2\sigma})  
\label{eq:bt2}\end{equation}
 for $ \pm \Im z\geq 0$ if $\tau <0$ and
 \begin{equation}
  \theta_{n} (z)= i  n^{-1/2} \sqrt{   \tau -   z n^{1-\sigma}}= 
   i \sqrt{\tau } n^{-1/2} - i \frac{   z}{2\sqrt{\tau}} n^{1/2-\sigma} + O(n^{3/2-2\sigma})
\label{eq:bt2+}\end{equation}
for all $z\in \Bbb C$ if  $\tau  > 0$.

 In  the case $\sigma <1$,  the term $ z n^{-\sigma}$  is dominating in \e{eq:ME1}. Moreover, for $\sigma\leq 2/3$, the phases $\theta_{n} (z)$ are  given by   formula  \e{eq:bt1-} more general than \e{eq:bt1}. The last circumstance is however inessential because the terms $t_{n}^l$ with $l>1$ in \e{eq:te} are negligible compared to $t_{n}$.  This yields an asymptotics 
   \begin{equation}
\theta_{n} (z)=\sqrt{ z} n^{-\sigma/2}
 \big(1+ O(n^{-\epsilon})\big), \q \epsilon>0.
 \label{eq:hm}\end{equation}


         
 In particular, these results imply the following assertion.

 \begin{proposition}\label{theta}
 Set 
 \begin{equation} 
\nu=
 \begin{cases}1/2 \q \mbox{if}\q\sigma \geq  1
\\ \sigma/2 \q \mbox{if}\q\sigma \leq  1
\end{cases}
\label{eq:nu}\end{equation}
and
 \begin{equation} 
\varkappa (z)= \begin{cases}\pm \sqrt{|\tau|} \q & \mbox{if}\q\sigma > 1, \,\tau<0, \, \pm \Im  z\geq 0
\\ 
i \sqrt{\tau} \q & \mbox{if}\q\sigma > 1,\,  \tau > 0, \, z\in {\Bbb C}
\\ 
 \sqrt{  z} \q &\mbox{if}\q\sigma < 1,   \, z\in  \clos \Pi_{0}, \, z\neq 0
 \\ 
 \sqrt{ z-\tau} \q &\mbox{if}\q\sigma = 1,    \, z\in  \tau + \clos \Pi_{0},  \, z\neq  \tau .
\end{cases}
\label{eq:wr}\end{equation}
Then
 \begin{equation} 
\theta_{n} (z)= \varkappa  (z)n^{-\nu} (1+ o (1)).
\label{eq:wra}\end{equation}            
              \end{proposition}

To pass to asymptotics of  sums \e{eq:sum}, we use 
     the Euler-Maclaurin formula
            \begin{equation}
 \sum_{m=1}^n F(m)= \int_{1}^n F(x) dx + \frac{F(n) +  F(1)}{2}+\int_{1}^n F'(x)\big( x-[x]- \frac{1}{2} \big) dx
\label{eq:E-M}\end{equation} 
where $[x]$ is the integer part of $x$.  This formula  is true for   arbitrary functions $F\in C^1$.

 Formula  \e{eq:E-M}   allows one to deduce an asymptotics as $n\to\infty$ of sum 
 \e{eq:sum} from that of the phases $\theta_{n}$.  For example, for $\sigma=1$, we apply \e{eq:E-M} to $F(x)= x^{-1/2}$ which yields 
            \begin{equation}
\varphi_{n}(  z )  = 2 \sqrt{-\tau +   z }\:  n^{1/2} +  C+ o(1)
\label{eq:xAf1}\end{equation} 
  with some constant $C$. The remainder $C+ o(1)$ here can be neglected in asymptotics \e{eq:A22G+} because the  Jost solutions are defined up to a constant factor.


  Next,  we  consider the case $\sigma\in (1,3/2)$. If $\tau< 0$, it follows from \e{eq:bt2} and the  Euler-Maclaurin formula \e{eq:E-M}  that  
      \begin{equation}
\varphi_{n}(  z )  = \pm 2 \sqrt{|\tau| n}  \pm  \frac{ z } {   \sqrt{ |\tau| } (3-2\sigma)}n^{3/2-\sigma}  + O (n^{5/2-2\sigma} )
\q\mbox{for} \pm \Im  z\geq 0.
\label{eq:Af}\end{equation} 
So,  up to error terms,  the functions $e^{i \varphi_{n}( z)}  $ where $z= \lambda + i\varepsilon $ contain oscillating 
\[
\exp\big(  
\pm 2 i  \sqrt{|\tau| n} \pm  \frac{i  \lambda} {   \sqrt{ |\tau| } (3-2\sigma)}n^{3/2-\sigma}  \big) 
\]
and  exponentially decaying\footnote{We say that a sequence $x_{n}$ tends to zero exponentially if  $x_{n}= O (e^{-n^a})$  for some $a>0$.}   
\[
\exp\big(  - \frac{|\varepsilon|} {   \sqrt{ |\tau| } (3-2\sigma)}n^{3/2-\sigma}  \big) 
\]
factors. Note that the strongly oscillating factor 
$\exp (  
\pm 2 i  \sqrt{|\tau| n } ) $ in the asymptotics of $f_{n} (z)$ as $n\to\infty$ does not depend on $z$. In the case $\tau >0$, we have
  \begin{equation}
  \varphi_{n}(  z )  = 2 i \sqrt{\tau  n}  -  \frac{ i z } {   \sqrt{ \tau } (3-2\sigma)}n^{3/2-\sigma} + O (n^{5/2-2\sigma} )\q \mbox{if}\q \sigma\in (1,3/2 ) .
\label{eq:Af+}\end{equation} 
Thus, the Jost solutions $f_{n}(  z )  $ contain an   exponentially decaying factor $e^{- 2  \sqrt{\tau  n}  }$  for all $z\in{\Bbb C}$. 

 Formulas \e{eq:Af} and \e{eq:Af+}    remain true also  for $\sigma=3/2$ if $(3-2\sigma)^{-1}n^{3/2-\sigma} $ is replaced by $\ln n$.  For example, for $\tau<0$, we have
   \begin{equation}
\varphi_{n}(  z )  = \pm 2 \sqrt{|\tau| n}  \pm  \frac{ z } {   \sqrt{ |\tau| }  } \ln n + C+ o (1 )
\q\mbox{for} \pm \Im  z\geq 0.
\label{eq:Aff}\end{equation}

 
 In the case $\sigma<1$,  asymptotics of the phases is given by relation \e{eq:hm}. 
      Therefore using formula \e{eq:E-M}, we find that
   \begin{equation}
\varphi_{n}(  z )  = 2 \sqrt{ z }(2-\sigma)^{-1} n^{1-\sigma/2}     + O (n^{\sigma/2}  ).
\label{eq:Af1}\end{equation} 
So, $e^{i\varphi_{n}(z)}$ exponentially decays if $ z\not\in [0,\infty)$ and   oscillates if $ z  = \lambda\pm i0$  for $\lambda>0$.
In the case $\sigma=1$, relation \e{eq:Af1} is true if $\sqrt{ z }$ is replaced by $\sqrt{-\tau+ z }$.

Note that explicit formulas for $\theta_{n} (z)$ allow one to find all power terms of asymptotic expansion of $\theta_{n} (z)$ as $n\to\infty$. In view of formula \e{eq:E-M} this yields all growing terms of the phases $\varphi_{n} (z)$ as $n\to\infty$.

 It follows from asymptotic formula  \e{eq:A22G+}  for the Jost solutions $f_{n} (z)$ and the results about the phases $\varphi_{n} (z)$  stated above that for all $\sigma\in (0,3/2)$ (for $\sigma>1$ it is also required that $\tau\neq 0$) and $\Im z\neq 0$, $f_{n} (z)$ tend to zero exponentially as $n \to \infty$.   In the critical case $\sigma=3/2$, the same is true  if $\tau >0$. If $\sigma=3/2$ and $\tau<0$, then relations   \e{eq:A22G+}, \e{eq:Aff} show that 
    \begin{equation}
f_{n} (\lambda+i\varepsilon) = (-\gamma)^n e^{ \pm 2 i\sqrt{|\tau| n} } n^{\pm i\gamma\lambda_{1}   } n^{-1/2-  \varepsilon_{1} }  \big(1+    o(1)   \big) \q\mbox{for} \pm   \gamma \varepsilon> 0,
\label{eq:Af1-}\end{equation} 
where $\lambda_{1} =\lambda /  \sqrt{ |\tau| } $, $\varepsilon_{1}=|\varepsilon| /  \sqrt{ |\tau| } $.

 In particular, we have

 \begin{proposition}\label{fl2}
 Under the assumptions of any of Theorems~\ref{GSS-x}, \ref{GSS-y}  or \ref{GSS-z}
 the inclusion 
   \begin{equation}
f_{n} (z) \in {\ell }^2  ({\Bbb Z}_{+}), \q   z \not\in \clos {\cal S},
\label{eq:Af1+}\end{equation} 
holds. In particular, \e{eq:Af1+}
is true  for $\Im z \neq 0$.           
              \end{proposition}
              
              Let us compare relation \e{eq:Af1-} with asymptotic formula (2.6) in \cite{crit}  for the singular case $\sigma>3/2$, $\tau <0$. The formula
               in \cite{crit}  is true for all $z\in{\Bbb C} $, the oscillating factor $e^{ \pm 2 i\sqrt{|\tau| n} } $ is the same as in   \e{eq:Af1-},  but the power of $n$ is $ 1/4 -\sigma/2$. This coincides with expression \e{eq:rho}, but  $1/4 -\sigma/2<-1/2$ for $\sigma>3/2$. In this case all solutions of equation \e{eq:Jy} are in ${\ell}^2 ({\Bbb Z}_{+} ) $ so that the deficiency indices of the operator $J_{\rm min}$ are $(1,1)$.
 
 Finally, we note that, on the absolutely continuous spectrum, formula  \e{eq:A22G+}  is consistent with a universal relation found in \cite{univ}.  Indeed, let   the assumptions of Theorems~\ref{GSS-x}, \ref{GSS-y} or \ref{GSS-z} be satisfied. 
 Using asymptotic formulas  \e{eq:xAf1}, \e{eq:Af}  or \e{eq:Af1}  and 
 calculating derivatives of the phases  $\varphi_{n} (\lambda\pm i0)$ in $\lambda$ we see  that,  with some constant factor $c_{\pm} (\lambda)$,
  \begin{equation}
d\varphi_{n} (\lambda\pm i0)/ d\lambda= c_{\pm} (\lambda)n^{\varsigma} (1+ o(1))
\label{eq:xAf1x}\end{equation} 
where  $\varsigma= 3/2-\sigma$ for $\sigma\in [1,3/2)$ and $\varsigma= 1-\sigma/2$ for $\sigma \leq 1$; if $\sigma =3/2$, then $n^{\varsigma}$ in  \e{eq:xAf1x} should be replaced by $\ln n$.  In view of definition \e{eq:rho}
 in all cases  the powers of $n$ in the amplitude and phase   in   formula    \e{eq:A22G+}  are linked by the equality
  \begin{equation}
2\rho+ \varsigma=1.
\label{eq:univ}\end{equation} 
This is one of the relations found in \cite{univ}; in the case $\sigma =3/2$, this relation reduces to the equality $\rho=1/2$.

For a comparison, we mention that, in the non-critical case $|\gamma|< 1$, we have $\rho=\sigma/2$ and $\varsigma=1-\sigma$ (see \cite{inf}) which is again consistent with equality \e{eq:univ}.

     \subsection{Exponentially growing solutions}
     

         For regular points $z\in{\Bbb C}$, the   solution $g_{n} (z)$ of equation   \e{eq:Jy}   linearly independent  with $f_{n} (z)$ is constructed by  formula \e{eq:GEg}. 
Using  the asymptotic formulas of Sect.~2.1 for the Jost solutions, we find a behavior  of $g_{n} (z)$ as $n\to\infty$.
  
     \begin{theorem}\label{P2+}
        Let one of the following three assumptions be satisfied:
        
       $1^0$ the conditions of Theorem~\ref{GSS-x} where either $\tau<0$, $\sigma<3/2$ and $\Im z \neq 0$ or 
        $\tau>0$ and  $z\in{\Bbb C}$ is arbitrary
       
       $2^0$ the  conditions of  Theorem~\ref{GSS-y} where either $\gamma =1$ and $ z\not\in [0,\infty)$ or 
        $\gamma = -1$ and $ z\not\in (-\infty, 0]$

    $3^0$   the conditions of  Theorem~\ref{GSS-z} where either $\gamma=1$ and $ z\not\in [\tau,\infty)$ or 
        $\gamma =-1 $ and $ z\not\in (-\infty, -\tau]$

Then   the  asymptotics of the solution $g_{n}(z) $ of equation   \e{eq:Jy}  is given by formula \e{eq:A2P3}.  In particular, 
    \begin{equation}
g_{n} (z) \notin {\ell }^2  ({\Bbb Z}_{+}) \q \mbox{if}\q z\not\in\clos{\cal S}.
 \label{eq:gL_2}\end{equation}
 \end{theorem}

 We emphasize that the definitions of the numbers $\rho$ and of the sequences $\varphi_{n} ( z )$ are different under  assumptions $1^0$, $2^0$ and $3^0$, but  asymptotic  formula  \e{eq:A2P3}  is true in all these cases. 
 
    In the critical case $\sigma=3/2$ (and $\tau <0$) the solution $g_{n} (z)$ of equation  \e{eq:Jy} behaves  as a power of $n$ as $n\to\infty$.

          \begin{proposition}\label{D1+} 
       If $\sigma=3/2$,  $\tau <0$ and $\pm   \gamma \varepsilon> 0$, then
       \begin{equation} 
g_{n} (\lambda+i\varepsilon) = (-\gamma)^n e^{ \pm 2 i\sqrt{|\tau| n} } n^{\pm i\gamma\lambda_{1} }  n^{-1/2+ \varepsilon_{1} }   \big(1+   o (  1 ) \big) 
\label{eq:gg+}\end{equation}
where $\lambda_{1} =\lambda /  \sqrt{ |\tau| } $, $\varepsilon_{1}=|\varepsilon| /  \sqrt{ |\tau| } $.
In particular,
 relation \e{eq:gL_2} is preserved.
      \end{proposition}
      
      Theorem~\ref{P2+}  and Proposition~\ref{D1+} will be proven in Sect.~6.1.
 

 All solutions  of equation   \e{eq:Jy}  and, in particular, the orthonormal  polynomials $P_{n}(z)$,   are linear combinations of the solutions  $f_{n} (z)$ and $g_{n} (z)$  for $z \not\in \clos{\cal S}$ or  of the solutions   $f_{n} (\lambda+ i0)$ and $f_{n} (\lambda-i0)$  for $z=\lambda\in{\cal S}$. Therefore the results stated above yield an asymptotics of $P_{n}(z)$ as $n\to\infty$. This is discussed in Sect.~6 -- see Theorem~\ref{P2}  and \ref{P1}.

   \subsection{Spectral results}
   
   First, we discuss the essential self-adjointness of    the minimal operator $J_{\rm min}$.  According to the limit point/circle theory this is equivalent to the existence of solutions of equation \e{eq:Jy}  where $\Im z\neq 0$ not belonging to  $\ell^2 ({\Bbb Z}_{+})$. Therefore the following result  is a direct consequence of Theorem~\ref{P2+} and Proposition~\ref{D1+}.
   
   \begin{proposition}\label{Self-Adj1}
 Let  assumptions  \e{eq:ASa},  \e{eq:ASb} with $| \gamma |=1$   and some $\sigma\in (0,3/2]$ be satisfied; for $\sigma>1$ we additionally suppose that $\tau \neq 0$.   
Then    the minimal operator $J_{\rm min}$ is essentially self-adjoint.   
         \end{proposition}

    Of course,  for $\sigma\leq 1$ one can refer to  the Carleman condition  \e{eq:Carl}, but   for $\sigma> 1$ the series in \e{eq:Carl}  is convergent.
   
   The case $\sigma>3/2$ was investigated in \cite{crit}, see Theorem~2.3.   According to Part~$2^0$ of this theorem, for $\tau>0$, 
   the  operator $J_{\rm min}$ remains essentially self-adjoint.  The results for the case  $\tau<0$ are more interesting. Combining   Proposition~\ref{Self-Adj1}  with 
   Part~$1^0$  of  Theorem~2.3 in \cite{crit}, we can state the following result.

 
     \begin{proposition}\label{Self-Adj}
 Suppose that assumptions  \e{eq:ASa},  \e{eq:ASb} with $| \gamma |=1$, $\tau<0$   and some $\sigma>0$ are satisfied. 
Then    the minimal operator $J_{\rm min}$ is essentially self-adjoint if and only if $\sigma\leq 3/2$.   
         \end{proposition}
         
         Note that Proposition~\ref{Self-Adj} does not contradict Theorem~2.1 of \cite{Domb} because the assumptions of \cite{Domb} correspond to the case
         $\tau=0$.
         
                  
         Below we  always suppose that $\sigma\leq 3/2$  and  denote by $J=\clos J_{\rm min}$  the closure of the  essentially self-adjoint operator $J_{\rm min}$.

   Spectral properties of Jacobi operators are determined by a behavior of solutions of equation  \e{eq:Jy} for real $z=\lambda$.  In particular, oscillating solutions correspond to the absolutely continuous spectrum.  On the contrary, for regular $\lambda$ or eigenvalues of $J$, one solution of   \e{eq:Jy}  exponentially decays and another one  exponentially grows. On the heuristic level, the results of Sect.~2.1 imply that the absolutely continuous spectrum of a  Jacobi operator $J$ consists of $\lambda$ where $-\tau n^{-1}+ \gamma \lambda n^{-\sigma} \geq 0$ (for large $n$).  On the contrary, the points $\lambda$   where  $-\tau n^{-1}+ \gamma \lambda n^{-\sigma}  < 0$ (again, for large $n$) are regular or,  eventually, are eigenvalues of $J$. This intuitive picture turns out to be correct.
   


  \begin{theorem}\label{S-Adj}

       Suppose that assumptions  \e{eq:ASa},  \e{eq:ASb} with $| \gamma |=1$   are satisfied. 
   
$1^0$
 Let  $\sigma\in (1,3/2 ]$. If $\tau <0$, then the spectrum of the operator $J$ is absolutely continuous and covers the whole real line.
If $\tau >0$, then the spectrum of the operator $J$ is discrete.

  $2^0$  Let  $\sigma\in (0,1)$. If $\gamma=1 $, then the  absolutely continuous spectrum of the operator $J$  coincides with the half-axis $[0, \infty)$ and its negative  spectrum of  $J$ is discrete. 
   If $\gamma=-1$, then the  absolutely continuous spectrum of the operator $J$  coincides with the half-axis $(- \infty, 0]$ and its positive  spectrum   is discrete.

    $3^0$  Let  $\sigma =1$. If $\gamma=1 $, then the  absolutely continuous spectrum of the operator $J$  coincides with the half-axis $[\tau, \infty)$ and its  spectrum  below the point $\tau$ is discrete.  If $\gamma=-1 $, then the  absolutely continuous spectrum of the operator $J$  coincides with the half-axis  $(-\infty,-\tau]$ and its  spectrum  above the point $-\tau$ is discrete. 
    
     \end{theorem}
    
    Parts $2^0$ and $3^0$ of Theorem~\ref{S-Adj}  can be considered as generalizations of the classical results  about  the Jacobi operators with the Laguerre coefficients \e{eq:Lag}.  We emphasize that, in the case $\sigma\leq 1$,  the are no conditions on the parameter $\tau$.
    The results of Part $1^0$ seem to be of a new nature.

 The results stated above apply to Jacobi operators with the coefficients $a_{n}$, $b_{n}$ growing as $n^\sigma$ where $\sigma$ is an arbitrary number in the interval $(0,3/2]$.  Together with the results of \cite{crit} where the case $\sigma>3/2$ was considered, they cover an arbitrary power growth of the Jacobi coefficients.

 Thus, our results show that,  in the critical case $|\gamma|=1$, there are two ``phase transitions": for $\sigma=1$ and  for $\sigma=3/2$. Indeed, the absolutely continuous spectrum of the Jacobi operator $J$ coincides with a half-axis for $\sigma\leq1$. In the case $\sigma\in (1,3/2]$, the spectrum 
 of  $J$ is either absolutely continuous   and  covers the whole real-axis for $\tau<0$ or it is discrete  for $\tau>0$.
   If $\sigma> 3/2$,  then the minimal Jacobi operator $J_{\rm min}$ has deficiency indices $(1,1)$ and
 the spectra of all its self-adjoint extensions are discrete.
 
 Our spectral results can be summarized in the following table where $\boldsymbol{\Sigma}_{\rm ac}$ and $\boldsymbol{\Sigma}_{\rm ess}$ are the absolutely continuous and  essential spectra of the operator $J$. For definiteness, we choose $\gamma=1$:

  \[
     \begin{matrix}
    \sigma\in (0,1 )& \Longrightarrow &\boldsymbol{\Sigma}_{\rm ac}=\boldsymbol{\Sigma}_{\rm ess}= [0,\infty)\\
    \sigma=1  &\Longrightarrow & \boldsymbol{\Sigma}_{\rm ac}=\boldsymbol{\Sigma}_{\rm ess}= [\tau,\infty)\\
      \sigma\in [1, 3/2], \; \tau <0 &\Longrightarrow &\boldsymbol{\Sigma}_{\rm ac} ={\Bbb R}\\
            \sigma\in [1, 3/2] , \; \tau >0 &\Longrightarrow &\boldsymbol{\Sigma}_{\rm ess} =\emptyset \\
                 \sigma > 3/2  &\Longrightarrow &\boldsymbol{\Sigma}_{\rm ess} =\emptyset 
 \end{matrix}
     \]
  

\section{Ansatz } 

 As usual,    we suppose that the  recurrence coefficients $a_{n}$, $b_{n}$ obey conditions \e{eq:ASa}, \e{eq:ASb} with $|\gamma| =1$. 
 We define the Ansatz $ {\cal A}_{n}= {\cal A}_{n}(z) $ by formula \e{eq:Gr1}  where the power $\rho$ and the phases $  \varphi_{n}= \varphi_{n}(\gamma z)$ will be found in this section.

   
   

 \subsection{Construction}

   Our goal  here is to  determine $\rho$ and $\varphi_{n}$ in such a way   that   remainder \e{eq:Grr} satisfies  condition \e{eq:rem}
for  
   \begin{equation}
   \d=1/2+\sigma
      \q \mbox{if} \q \sigma > 1 \q \mbox{and some}\footnote{The precise value of $\d=\d(\sigma)$ for  $\sigma \in(2/3,1)$ is indicated in Propositions~\ref{AA+-}. For $\sigma\leq 2/3$, it can be deduced from the proof of Proposition~\ref{AA-}, but we do not need it.}
        \q \d > 1+\sigma/2  \q \mbox{if} \q \sigma<1.
\label{eq:R1}\end{equation}
If $\sigma=1$, then $\d=2$, so that the estimate of the remainder is more precise in this particular case.
  We emphasize   that   estimate \e{eq:rem} with $\d>1$   used  in the non-critical case $| \gamma | \neq 1$ in \cite{inf}  is not sufficient now.



    Put
    \begin{equation}
{\cal B}_{n}=   \frac{{\cal A}_{n+1}}{{\cal A}_{n}}.
\label{eq:R2}\end{equation}
Then  expression \e{eq:Grr}  for the remainder can be rewritten as
 \begin{equation}
{ \bf r}_{n} (z) = \sqrt{ \frac {a_{n-1}} {a_{n}}} {\cal B}_{n-1}^{-1}
 +
  \sqrt{ \frac {a_{n}} {a_{n-1}}}  {\cal B}_{n}
   + 2 \gamma_{n} -  \frac{z}{\sqrt{ a_{n-1}a_{n}}} .
\label{eq:R3}\end{equation}

Assumption   \e{eq:ASa} on $a_{n}$ implies  that
\begin{equation}
  \sqrt{ \frac {a_{n}} {a_{n-1}}}  = (n+1)^{\sigma/2} n^{-\sigma/2}(1 + O (n^{-2}) )=1+ \frac{\sigma/2} {n}+ O (n^{-2}) 
\label{eq:BX}\end{equation}
and
\[
(a_{n} a_{n-1})^{-1/2}=   n^{-\sigma} \big(1 + O (n^{-1}) \big).            
\]
Using also assumption \e{eq:ASb} on $b_{n}$, we see that  sequence  \e{eq:Gr}  satisfies a relation
\begin{equation}
 \gamma_{n}=  1+  (\tau/2 ) n^{-1}+ O (n^{-2}) 
\label{eq:BX1}\end{equation}
where $\tau$ is defined by equality
\e{eq:BX3}.

We seek $ {\cal A}_{n}$ in form \e{eq:Gr1} where the phases $\varphi_{n}$ are defined as sums   \e{eq:sum}.  The power $\rho$ and the differences  
 \[
 \theta_{n}  = \varphi_{n+1}- \varphi_{n}
\]
will be determined by condition  \e{eq:rem}.  The  sequences $ \theta_{n} $ constructed below  tend to zero as $n\to\infty$ and satisfy  condition  \e{eq:the}.
It follows from \e{eq:Gr1} and \e{eq:R2} that
 \begin{equation}
{\cal B}_{n}   = -  (n+1)^{-\rho}  n^\rho
e^{ i\theta_{n}  } = - (1 -\rho  n^{-1} + O (n^{-2}) )
e^{ i\theta_{n}  } .
\label{eq:Grq}\end{equation}

According to  relations   \e{eq:BX} -- \e{eq:BX1}  and \e{eq:Grq} the following intermediary assertion is a direct consequence of expression \e{eq:R3}.

\begin{lemma}\label{RK}
Relative remainder \e{eq:Grr} admits a representation
  \begin{equation}
 {\bf r}_{n}= -   \big(1 - (\nu/2) n^{-1}  \big)  e^{- i\theta_{n-1}  } 
-   \big(1 +  (\nu/2)  n^{-1}  \big)    e^{ i\theta_n  } + 2 + \tau n^{-1} - z n^{-\sigma} + O (n^{-  \d})
\label{eq:R3a}\end{equation}
where  
  \begin{equation}
\nu= \sigma - 2 \rho
\label{eq:R3aa}\end{equation}
 and $\d =\min\{2, 1+\sigma\}$.
 \end{lemma}
 
 Note that in view of \e{eq:rho} expressions  \e{eq:nu} and \e{eq:R3aa}  for
 $\nu$ are equivalent.
  
 For all $\sigma\in (2/3, 3/2]$, the phases  $\theta_{n} $  are defined by  the same formulas \e{eq:ME1} and \e{eq:bt1}, that is,
  \begin{equation}
\theta_{n} =  \theta_{n}  (z)   =\sqrt{ -\tau n^{-1} +  z n^{-\sigma} } ,\q \Im \theta_{n} (z) \geq 0,
\label{eq:t=th}\end{equation}
although the estimates of the remainder $ {\bf r}_{n}$ are rather different in the cases  $\sigma > 1$, $\sigma=1$ and $\sigma < 1$. For $\sigma\leq 2/3$, expression \e{eq:t=th}   requires some corrections.

 \subsection{The case $\sigma > 1$ }
 For such $\sigma$, we  suppose that $\tau\neq 0$.   We treat the cases $\tau <0$ and $\tau >0$ parallelly  putting  $\sqrt{- \tau } >0$ if $  \tau < 0$  and  $\sqrt{ -\tau } = i\sqrt{| \tau |}$ if $  \tau >0$.

It follows from definition \e{eq:t=th} that
$\theta_{n}=O(n^{-{1/2}})$,
whence
 \begin{equation}
e^{i\theta_n} =\sum_{k=0}^3 \frac{i^k}{k!}
 \theta_{n}^k+ O(n^{-2}).
\label{eq:R5}\end{equation} 
 Substituting \e{eq:R5} into representation \e{eq:R3a}, we see that
  \begin{equation}
{\bf r}_{n}=  \sum_{k=0}^3   r_{n}^{(k)}
 + O (n^{- 2})
 \label{eq:RX}\end{equation} 
 where
   \begin{align}
r_{n}^{(0)}=& -   \big(1 + (\nu /2) n^{-1}  \big)   
-   \big(1 -  (\nu /2)  n^{-1}  \big)    + 2 + \tau n^{-1} - z n^{-\sigma} 
 \nonumber\\
 = &\tau n^{-1} -  z n^{-\sigma} =-t_{n},
 \label{eq:RX1}\end{align}
   by definition  \e{eq:ME1}, 
 and
   \begin{align}
r_{n}^{(1)}= & i  \big(1 -  (\nu /2) n^{-1}  \big)   \theta_{n-1}
- i  \big(1 +   (\nu /2) n^{-1}  \big)   \theta_{n},
 \label{eq:RX1A}\\
2 r_{n}^{(2)}=  & \big(1 -  (\nu /2)  n^{-1}  \big)   \theta_{n-1}^2
 +  \big(1 +   (\nu /2) n^{-1}  \big)   \theta_{n}^2,
 \label{eq:RX1B}\\
6 r_{n}^{(3)}= & - i  \big(1 -  (\nu /2) n^{-1}  \big)   \theta_{n-1}^3
+ i  \big(1  +  (\nu /2)   n^{-1}  \big)   \theta_{n}^3 .
 \label{eq:RX1C}\end{align}

   Since $ \theta_{n}^2=t_{n}$, it follows from \e{eq:RX1B}   that
      \begin{align}
2 r_{n}^{(2)}=&    \big(1 -  (\nu /2) n^{-1}  \big)   t_{n-1}
 +  \big(1  + (\nu /2)n^{-1}  \big)   t_{n} 
 \nonumber\\
 = & 2 t_{n} +  \big(1 - (\nu /2) n^{-1}  \big)  (t_{n-1} -t_{n}).
 \label{eq:RX1B2}\end{align} 
   Comparing this equality with \e{eq:RX1}, we find that 
     \begin{equation}
r_{n}^{(0)} + r_{n}^{(2)}= 2^{-1}  \big(1 - (\nu /2) n^{-1}  \big)  (t_{n-1} -t_{n})=  O(n^{-2}) .
  \label{eq:ME7xx}\end{equation}


 The power $\rho$ in \e{eq:Gr1} is determined by
   linear term \e{eq:RX1A}  which we write as
      \begin{equation}
   r^{(1)}_{n} = i(  \theta_{n-1} -\theta_{n} ) - i (\nu /2) n^{-1} ( \theta_{n}  + \theta_{n-1}) .
 \label{eq:ME7x}\end{equation} 
 Let us distinguish   the leading term in \e{eq:t=th}  setting
   \begin{equation}
\theta_{n} =  \sqrt{ -\tau n^{-1} } + \ti\theta_{n}
\label{eq:t=th1}\end{equation}
where
 \begin{equation}
\ti\theta_{n} =  \sqrt{ -\tau n^{-1} + z n^{-\sigma}} - \sqrt{ -\tau n^{-1}  }=\frac{zn^{1/2-\sigma}}{  \sqrt{ -\tau  + z n^{1-\sigma}} + \sqrt{ -\tau   }}= O(n^{1/2-\sigma}).
\label{eq:t=th2}\end{equation}
Let us substitute \e{eq:t=th1} into \e{eq:ME7x} and observe that
\[
(\sqrt{  (n-1)^{-1} }-\sqrt{  n^{-1} })-  (\nu/2) n^{-1}(\sqrt{  (n-1)^{-1} }+\sqrt{  n^{-1} })=(2^{-1}-\nu)n^{-3/2}+ O (n^{-5/2}).
\]
According to \e{eq:t=th2} we have 
 \begin{equation}
\ti\theta_{n} - \ti\theta_{n-1} = O(n^{-\sigma-1/2}).
\label{eq:hold}\end{equation}
Thus, it follows from \e{eq:ME7x} that
    \[
  r^{(1)}_{n} = i   \sqrt{-\tau} (  2^{-1} - \nu) n^{-3/2}
+ O(n^{-\sigma-1/2}).
\]
  The  coefficient of $n^{-3/ 2} $ here is zero if $\nu=1/2$ which, by \e{eq:R3aa},  yields
   $ \rho= \sigma/2 - 1/4$; in this case 
$ r^{(1)}_{n} = O (n^{- \sigma -1 /2} )$.  



   
   It remains to consider the term    $r_{n}^{(3)}$. In  view of \e{eq:RX1C}, it equals 
      \begin{equation}
6 r_{n}^{(3)}=   i  ( \theta_{n}^3- \theta_{n-1}^3)  + i  ( \nu /2) n^{-1} ( \theta_{n}^3+ \theta_{n-1}^3)  .
 \label{eq:RX1Ct}\end{equation}  
   Observe that 
     \begin{equation}
  \theta_{n}^3 -    \theta^3 _{n-1}=  (\theta_{n } -    \theta_{n-1}) (\theta_{n }^2+ \theta_{n }   \theta_{n-1}  +   \theta_{n-1}^2).
  \label{eq:r3-}\end{equation}  
It follows from relations \e{eq:t=th1}  and \e{eq:hold} that 
  the first factor here is $O(n^{-3/2})$. The second factor is $O(n^{ - 1})$ because $\theta_{n}=O(n^{ - 1/2})$.  Therefore expression \e{eq:r3-} is 
  $O(n^{-5/2})$. Obviously, 
  the second term in the right-hand side of \e{eq:RX1Ct}  satisfies the same estimate.


  Let us state the result obtained. 
 
 \begin{proposition}\label{AA+}
 Let the assumptions of Theorem~\ref{GSS-x}  be satisfied, and let   the phases $\theta_{n} (z)$ be given by formula  \e{eq:t=th}.
Define the Ansatz ${\cal A}_{n} (z)$  by formula \e{eq:Gr1} where $ \rho= \sigma/2 - 1/4 $. Then   remainder \e{eq:Grr}  satisfies estimate \e{eq:rem} where $\d=\sigma+1/2$.
 \end{proposition}

  \subsection{The intermediary case $\sigma  = 1$ }
  
  The results of this subsection are a particular case of Proposition~\ref{AA+}, but the construction of the phases is  now simpler:
     \begin{equation}
t_{n}  = ( z-\tau) n^{-1}  \q \mbox{and} \q \theta_{n} = \sqrt{ z-\tau} n^{-1/2} .
 \label{eq:La}\end{equation}
The estimate of
   the remainder ${\bf r}_{n} $  is also simpler and more precise than in the general case. Indeed,   according to  \e{eq:La}, we now have
   \[
 \theta_{n-1}- \theta_{n} = 2^{-1}\sqrt{ z-\tau} n^{-3/2} + O (n^{-5/2}).
 \]
 Therefore,
 it follows from \e{eq:RX1A} where $\nu =1/2$ that  $ r_{n}^{(1)}= O (n^{-5/2})$. The same estimate for   $ r_{n}^{(3)}$ is a direct consequence of  \e{eq:r3-}. Estimate \e{eq:ME7xx}  remains of course true.
  Thus, using equality \e{eq:RX} we can state the limit case of Proposition~\ref{AA+}.
  
   \begin{proposition}\label{AA+1}
 Let the assumptions of Theorem~\ref{GSS-z}  be satisfied, and let   the phases $\theta_{n}(z)$ be given by formula  \e{eq:La}.
Define the Ansatz ${\cal A}_{n}(z)$  by formula \e{eq:Gr1} where $ \rho=  1/4 $. 
Then  remainder \e{eq:Grr}  satisfies estimate \e{eq:rem} where $\d=2$.
 \end{proposition}

 \subsection{The case $\sigma \in (2/3,1)$ }
 
 We again define the phases $\theta_{n}  $ by formula \e{eq:t=th}, but now
the term $ z n^{-\sigma}$  is dominating    so that, instead of  \e{eq:t=th1},  \e{eq:t=th2}, we have a relation   
 \begin{equation}
\theta_{n}  =\sqrt{ t_{n}   }  =\sqrt{   z  } n^{-\sigma /2} \big( 1+ O ( n^{\sigma -1})\big).
\label{eq:tt}\end{equation}
Therefore  
the scheme exposed  in  Sect.~3.2 for the case $\sigma>1$ requires some modifications.

It again suffices to keep $4$ terms in expansion of $e^{i\theta_{n}}$, but the remainders in formulas  \e{eq:R5} and \e{eq:RX} are now $O (n^{-2\sigma})$. 
Estimates of $r_{n}^{(k)}$ where $k=0,1,2,3$ are the same as in Sect.~3.2  if the roles of the terms $-\tau n^{-1}$ and 
$ z n^{-\sigma}$ are interchanged. 
Relations \e{eq:RX1} and \e{eq:RX1B2}  are preserved, but the remainder $O (n^{-2 })$ in \e{eq:ME7xx}  is replaced by 
$O (n^{-1-\sigma})$.  It directly follows from definition \e{eq:ME1}  that
   \begin{equation}
  t_{n-1} -t_{n}= z \sigma n^{-1-\sigma} \big(1+ O(n^{-1+\sigma})\big).
 \label{eq:hold3}\end{equation} 
 Similarly to \e{eq:hold}, it follows from 
 \e{eq:tt}, \e{eq:hold3} that
   \begin{equation}
 \theta_{n-1} -\theta_{n}= 2^{-1}\sqrt{ z} \sigma n^{-1-\sigma/2}\big(1+ O(n^{-1+\sigma})\big).
 \label{eq:hold1}\end{equation} 
Therefore expression \e{eq:ME7x} equals
    \begin{equation}
  r^{(1)}_{n} = i    \sqrt{ z} (\sigma/2- \nu) n^{-1-\sigma/2}
+ O(n^{-2+\sigma/2}).
 \label{eq:rem1}\end{equation} 
  The coefficient at $n^{-1- \sigma/2} $ is zero if $\nu= \sigma/2$ which yields $2\rho=\sigma-\nu=\sigma/2$; in this case 
$ r^{(1)}_{n} = O(n^{-2+\sigma/2})$. Putting together equality  \e{eq:ME7xx}  and estimate  \e{eq:hold3} , we see that
$r_{n}^{(0)} + r_{n}^{(2)} =  O(n^{-1-\sigma}) $  which is $O(n^{-2\sigma}) $ because $\sigma<1$.
According to \e{eq:tt} and \e{eq:hold1} expression \e{eq:r3-} is estimated by $C n^{-1-3\sigma /2}$. In view of \e{eq:RX1Ct}  the same bound is true for $r_{n}^{(3)}$.

    Thus, we arrive at the following result. 
 
  \begin{proposition}\label{AA+-}
 Let the assumptions of Theorem~\ref{GSS-y}  be satisfied  with $\sigma\in (2/3, 1)$, and let  the phases $\theta_{n} (z)$ be given by formula  \e{eq:t=th}. 
 Define the Ansatz ${\cal A}_{n} (z)$  by formula  \e{eq:Gr1}  where $\rho=\sigma/4$.  Then  remainder \e{eq:Grr}  satisfies estimate \e{eq:rem}  with 
  $\d =\min\{2\sigma,  2- \sigma/2 \}> 1+\sigma/2$.
 \end{proposition}
 
We emphasize that for all  $\sigma\in (2/3, 3/2]$ the phases $\theta_{n}$ are given by the same formula \e{eq:t=th}. However asymptotics of $\theta_{n} $ are different for $\sigma>1$ and for $\sigma<1$ -- 
cf. \e{eq:t=th1}, \e{eq:t=th2}   with \e{eq:tt}.

 \subsection{The case $\sigma \leq 2/3$. Eikonal equation. }
 
 The leading term of the asymptotics of the phases $\theta_{n} $ is again given by formula \e{eq:tt}, but, additionally,  lower order terms appear.
Now,
we need to keep more terms in expansion \e{eq:R5} setting
 \begin{equation}
e^{i\theta_n} =\sum_{k=0}^{K} \frac{i^k}{k!}
 \theta_{n}^k+ O(n^{-  (K+1) \sigma /2})  .
\label{eq:R5+}\end{equation} 
Substituting \e{eq:R5+} into representation  \e{eq:R3a}, we see that
  \begin{equation}
 {\bf r}_{n}=  \sum_{k=0}^{K}  r_{n}^{(k)}
 + O(n^{-  (K+1) \sigma /2}) ,
 \label{eq:RX+}\end{equation} 
 where   $ r_{n}^{(0)}$ are again given by equality \e{eq:RX1}  and  
     \begin{align}
- i^k   k!  r_{n}^{(k)}=&   \big(1 -  ( \nu /2) n^{-1} \big)  \theta_{n-1}^k  +  (-1)
 ^k  \big(1 +  ( \nu /2)  n^{-1}  \big)  \theta_{n}^k
\nonumber\\
=&\theta_{n-1}^k + (-1)^k \theta_{n}^k -  ( \nu /2) n^{-1} \big(\theta_{n-1}^k - (-1)^k \theta_{n}^k  \big),
 \q k\geq 1 .
 \label{eq:RX2}\end{align}
 Of course, for $k=1,2,3$, this expression coincides with  \e{eq:RX1A},  \e{eq:RX1B},  \e{eq:RX1C}, respectively.
 It is convenient to choose an even $K=2L$ with a sufficiently large $L$.  We suppose that
   \begin{equation}
    (L+1/2)\sigma >1. 
   \label{eq:LL}\end{equation}

 Let us distinguish   the terms corresponding to $k=0$ and $k=1$ in  sum \e{eq:RX+} and then split it into the sums over even and odd $k$: 
    \[
   {\bf r}_{n}=  r^{(0)}_{n}  + r^{(1)}_{n}  + {\bf r}^{(ev)}_{n} + {\bf r}^{(odd)}_{n} + O (n^{-(L+1/2)\sigma}), 
 \]
 where 
 $ r^{(0)}_{n}  $, $r^{(1)}_{n} $ are given by formulas \e{eq:RX1},  \e{eq:ME7x} and
    \begin{equation}
  {\bf r}^{(ev)}_{n}= \sum_{l=1}^{ L }   r_{n}^{(2l)}, \q   {\bf r}^{(odd)}_{n}= \sum_{l=1}^{L -1}   r_{n}^{(2l+1)}.
 \label{eq:RXe}\end{equation} 
 
 To satisfy estimate \e{eq:rem} with a suitable $\d$, we now have to take the even terms $r_{n}^{(2l)}$ for all $l\leq L$ into account. The odd terms $r_{n}^{(2l+1)}$ turn out to be negligible. To be precise,
 we define the  phases $\theta_{n}$   by formula \e{eq:bt1-}  where $T_{n}$ is sum \e{eq:te}.  The coefficients $p_{l}$ will be found from the relation
     \begin{equation}
       r^{(0)}_{n}  +   {\bf r}^{(ev)}_{n} = O (n^{-1-\sigma}) 
        \label{eq:rr}\end{equation} 
        generalizing \e{eq:ME7xx}. To satisfy this relation, we use that the differences  between $\theta_{n}$ and $\theta_{n-1}$ in  the expression
           \[
(-1)^{l+1}  (2l)!  r_{n}^{(2l)} 
=\theta_{n-1}^{2l} + \theta_{n}^{2l} 
 -  ( \nu /2) n^{-1} \big(\theta_{n-1}^{2l} -  \theta_{n}^{2l} \big)
 \] 
    (it is a particular case of  \e{eq:RX2})   can be neglected.  Thus, we set
         \[
\Theta_{n}= 2\sum_{l=1}^{L} \frac{(-1)^{l+1}} {(2l)!} \theta_{n}^{2l}.
  \]
Since $ r^{(0)}_{n}=-t_{n}$, we find  that
    \begin{equation}
  r^{(0)}_{n}  +   {\bf r}^{(ev)}_{n} = (- t_{n} + \Theta_{n} ) + (1- (\nu/2) n^{-1} ) \sum_{l=1}^L (-1)^l (2l)!^{-1}(\theta_{n}^{2l} -\theta_{n-1}^{2l}).
     \label{eq:rr1}\end{equation} 
     As we  will see the sum here is negligible, and hence
      we can replace \e{eq:rr} by the (approximate) eikonal equation
         \begin{equation}
\Theta_{n}= t_{n}  +O (n^{-1-\sigma}). 
 \label{eq:MEm}\end{equation}
 
 Our goal is to   solve this equation with respect to $\theta_{n}^{2} $.
 Note that  $\Theta_{n}=\theta_{n}^2 $ if $L=1$ so that \e{eq:MEm} again yields expression $\theta_{n}^2= t_{n} $.
  The following elementary assertion shows that  equation \e{eq:MEm} can be efficiently solved for all $L\geq 1$. It is convenient to consider this problem in a somewhat more general setting.
  Denote by ${\cal P}$ the set of all polynomials (of the variable $t$), and let  ${\cal P}_{L}= t^{L+1} {\cal P}$, that is, ${\cal P}_{L}  \subset {\cal P}$ 
  consists of polynomials with zero coefficients at powers $t^k$  for all $k=0,1,\ldots, L$.

\begin{lemma}\label{ME}
Let $L\geq 2$ and $a_{2}, \ldots, a_{L}$ be arbitrary given  numbers. Then there exists a polynomial
  \[
P_{L} (t)=\sum_{l=2}^L p_{l} t^l 
 \]
  such that  the polynomial
   \begin{equation}
Q_{L}(t ) := P_{L} (t) +\sum_{k=2}^L a_{k} (P_{L} (t) +t)^k \in {\cal P}_{L}.
 \label{eq:al}\end{equation}
   \end{lemma}

  \begin{pf}   
  For arbitrary $ p_{2}, \ldots, p_{L}$, the polynomial $Q_{L}(t)$ defined by \e{eq:al} has degree $L^2$ and it does not contain  terms with zero and first powers of $t$.  We have to choose the numbers $ p_{2}, \ldots, p_{L}$ in such a way that   the coefficients of $Q_{L}(t)$ at $t^l$  are zeros for all $l=2,\ldots, L$.   This assertion is obvious for $L=2$ because
  \[
  Q_{2} (t)= P_{2}  (t) +a_{2} (P_{2}  (t) +t)^2
 = (p_{2}+a_{2})t^2 +2 a_{2} p_{2}t^3+a_{2}p_2^2 t^4,
  \]
  and, so,  $Q_{2}  (t) \in {\cal P}_{2}$ if $p_{2}=-a_{2}$.
  
  Let us pass to the general case. Suppose that \e{eq:al} is satisfied. Then there exists a number $q_{L+1}$ such that
          \begin{equation}
Q_{L}(t) -q_{L+1}t^{L+1}\in {\cal P}_{L+1}.
  \label{eq:QL}\end{equation}  
  We will find a number $ p_{L+1} $ such that
  the  polynomial
     \begin{equation}
  P_{L+1} (t)=P_{L} (t) + p_{L+1}   t^{L+1}
  \label{eq:al2}\end{equation}
  satisfies   \e{eq:al}   for $L+1$, that is, 
       \begin{equation}
Q_{L+1}(t ) : = P_{L+1} (t) +\sum_{k=2}^{L+1}a_{k} (P_{L+1} (t) +t)^k \in {\cal P}_{L+1}.
 \label{eq:al+1}\end{equation}  
 Let us calculate the polynomial $Q_{L+1}(t)$ neglecting  terms  in ${\cal P}_{L+1}$.  First, we observe
   that,   for all $k=2,\ldots, L,L+1$, the difference
     \[
( P_{L+1}  (t) +t)^k-  ( P_{L}  (t) +t)^k = \sum_{n=1}^k  \binom{k}n p_{L+1}^n  t^{(L+1)n} ( P_{L}  (t) +t)^{k-n} \in {\cal P}_{L+1}.
 \]
 Using also \e{eq:al2}, we see that, up to terms in ${\cal P}_{L+1}$, polynomial \e{eq:al+1} equals
   \[
Q_{L+1}(t)= P_{L} (t) + p_{L+1}   t^{L+1}+   \sum_{k=2}^{L} a_{k} (P_{L} (t) +t)^k  + a_{L+1} (P_{L} (t) +t)^{L+1}
 \]
  whence,  by assumption \e{eq:al},  
   \[
Q_{L+1}(t)= Q_{L}(t) + (p_{L+1}  + a_{L+1} ) t^{L+1} \in {\cal P}_{L+1}.
 \]
 It follows from \e{eq:QL}  that  this relation   is equivalent to
   \[
Q_{L+1}(t)- (p_{L+1} +q_{L+1}  + a_{L+1}) t^{L+1}\in {\cal P}_{L+1}.
 \]
 Thus,  inclusion  $Q_{L+1}(t) \in {\cal P}_{L+1}$  is true  if $p_{L+1} =-a_{L+1} - q_{L+1} $. This proves \e{eq:al} for $L+1$.
   \end{pf}

   Note particular cases
    \[
   p_{2}=-a_{2}, \q p_{3}=   2a_{2}^2 -a_{3}.
     \]

  Let us come back to relation \e{eq:MEm}.  Let us use Lemma~\ref{ME} with the coefficients $a_{l}=2 (-1)^{l+1} /(2l)! $, $t=t_{n}$ defined by equality \e{eq:ME1}, and let $p_{l}$ be  the coefficients constructed in this lemma.  It follows from  equality \e{eq:al}  that   the phases
     \begin{equation}
\theta_{n}^2 = t_{n} + \sum_{l=2}^L p_l t_{n}^{l} =: T_{n}
 \label{eq:ME2}\end{equation}
 satisfy, for some coefficients $q_{l}$, the equation
  \begin{equation}
2\sum_{l=1}^{L} \frac{(-1)^{l+1}} {(2l)!} \theta_{n}^{2l}   - t_{n} = \sum_{l=L+1}^{L^2} q_{l}t_{n}^l= O(t_{n}^{-L -1}).
 \label{eq:MDB}\end{equation}
  Since $t_{n}= O(n^{-\sigma})$, the right-hand side here  is $  O(n^{- (L+1)\sigma})$  which is $  O(n^{-  \d })$ with $\d> 1+ \sigma/2$
  if   condition  \e{eq:LL}  is satisfied.

 The definition of the phases by formula \e{eq:ME2} coincides of course with their definition by relations \e{eq:te},  \e{eq:bt1-}.
The asymptotics of $\theta_{n} $ as $n\to\infty$ is given by formula
\e{eq:hm}
 generalizing \e{eq:tt}. Next, we estimate the differences
  \[
 \theta_{n-1} -\theta_{n}=\frac{ T_{n-1} -T_{n}}{ \theta_{n-1} +\theta_{n}} .
 \]
According to \e{eq:te} we have
   \[
T_{n-1} -T_{n}= t_{n-1} -t_{n}+ \sum_{l=2}^{L} p_{l} ( t_{n-1}^l -t_{n}^l)
 \]
 so that it satisfies the same relation \e{eq:hold3} as $t_{n}$:
  \[
  T_{n-1} -T_{n}= z \sigma n^{-1-\sigma} \big(1+ O(n^{-1+\sigma})\big).
 \]
  Combining this relation  with \e{eq:hm}, we see that 
    \begin{equation}
 \theta_{n-1} -\theta_{n}=  2^{-1}\sqrt{ z} \sigma n^{-1-\sigma/2}\big(1+ O(n^{-\epsilon})\big)
 \label{eq:hold11}\end{equation} 
 for some $\epsilon>0$
 (compared with \e{eq:hold1} only the estimate of the remainder is changed).
 
It easily follows from \e{eq:hm} and \e{eq:hold11} that
   \begin{equation} 
|\theta_{n-1}^k-\theta_{n}^k | \leq C_{k} n^{-1 -k\sigma/2}  
 \label{eq:MB2-}\end{equation}
 for all $k=1, 2, \ldots $.


    Let us come back to Ansatz \e{eq:Gr1}. Similarly to Sect.~3.4,                                                                                                                                                         
     the power $\rho$ in \e{eq:Gr1} is determined by
  the linear term $r^{(1)}_{n}$ given by equality
    \e{eq:ME7x}. It again satisfies relation \e{eq:rem1} (with the remainder $O(n^{-2+\sigma/2})$ replaced by $O(n^{-\d})$ for some  $\d>1+\sigma/2$).    The coefficient at $n^{-1- \sigma/2} $ is zero if $\nu = \sigma/2$ which yields $\rho=\sigma/4$; in this case 
$
 r^{(1)}_{n} = O(n^{-\d}).
$


 Given  inequalities \e{eq:hm} and \e{eq:MB2-}, we can estimate the remainder ${\bf r}_{n}$ essentially similarly to Proposition~\ref{AA+-}.
   The only differences are that estimates of the remainders are slightly  weaker and that we have to take into account higher powers of $\theta_{n}$.
   First, we consider term  \e{eq:rr1}  with even powers of $\theta_{n}$. Both  the first term $- t_{n} + \Theta_{n} $ and the sum on the right are  $O(n^{-1-\sigma}) $
by virtue of relations \e{eq:MDB}  and  \e{eq:MB2-}, respectively. 
        The term
 ${\bf r}^{(odd)}_{n}$ is also negligible. Indeed, according to  \e{eq:RX2}  and \e{eq:RXe}  it equals
 \[
{\bf r}^{(odd)}_{n}=i\sum_{l=1}^{L-1} \frac{(-1)^l} {(2l+1)!}\Big((\theta_{n-1}^{2l+1}- \theta_{n}^{2l+1}) -(\nu/2)  n^{-1} (\theta_{n-1}^{2l+1}+ \theta_{n}^{2l+1})\Big).
\]
 Relations \e{eq:hm}  and  \e{eq:MB2-} allow us to estimate all terms here by $n^{-1-3\sigma/2}$.
 

Thus, we arrive at the following assertion generalizing Proposition~\ref{AA+-}.
  
   \begin{proposition}\label{AA-}
 Let the assumptions of Theorem~\ref{GSS-y}  be satisfied, and let   the phases $\theta_{n}(z)$ be given by  formulas  \e{eq:te}, \e{eq:bt1-} with $(L+1/2)\sigma >1$ and the coefficients $p_{2},\ldots, p_{L}$ constructed in Lemma~\ref{ME}.
Let the Ansatz ${\cal A}_{n}(z)$ be defined by formula \e{eq:Gr1} where $ \rho= \sigma/4 $. Then   remainder \e{eq:Grr}  satisfies estimate \e{eq:rem} with some $\d > 1+ \sigma/2 $.
 \end{proposition}
   
  \begin{remark}\label{the}
  In all estimates, we suppose that $ z\in    \clos \Pi_{0}$, $0<r\leq |z|\leq R<\infty$ for some $r$ and $R$ and $n\geq N = N (r,R)$. Then it follows from equality
  \e{eq:te}  that $\pm \Im T_{n} (z)>0$ as long as $\pm \Im t_{n} (z)>0$, that is, $\pm \Im  z >0$.
   Therefore 
  $T_{n} (z)\in    \clos \Pi_{0}$, and hence   condition \e{eq:the} is satisfied.
   \end{remark}
   
 Note two particular cases. If $\sigma>2/3$, then
   we can take $L=1$; this is the case considered in Proposition~\ref{AA+-}. If $\sigma>2/5$, then 
  $L=2$  so that  the formula for $\theta_{n}$ contains only one additional (compared with  \e{eq:bt1}) term:
\[
\theta_{n}=\sqrt{t_{n}  + t_{n}^2 /6}.
\]

 We, finally, note that constructions of  Ans\"atzen  were important steps also in the papers \cite{Jan-Nab-Sh, Na-Si}. However, the form of the  Ansatz ${\cal A}_{n}(z)$ suggested in this section is   different  from  \cite{Jan-Nab-Sh, Na-Si}; in particular, the phases $\varphi_{n} (z)$ in \e{eq:Gr1} are simplest in the case $\sigma>2/3$ while this case was excluded in  \cite{Jan-Nab-Sh}.   They are also different from \cite{Na-Si} -- see Remark~\ref{Na-Si}.

      \section{Difference and Volterra equations}

 Here we reduce a construction of the Jost solutions $f_{n}  (z)$ of the Jacobi equation \e{eq:Jy} to a Volterra ``integral" equation which is then solved by iterations.   In this section,  we do not make any specific assumptions about  the recurrence coefficients $a_{n}$, $b_{n}$ and  the Ansatz ${\cal A}_{n}(z)$ except of course that ${\cal A}_{n} (z)  \neq 0$; for definiteness, we set ${\cal A}_{-1}=1$.  We present a  general scheme of investigation and then,   in Section~5, apply it to Jacobi operators with  coefficients $a_{n} $ and $b_{n}$ satisfying conditions \e{eq:ASa} and \e{eq:ASb}  with $|\gamma|=1$.

   \subsection{Multiplicative change of variables}   
   
    For a construction of $f_{n}  (z)$, we will reformulate the problem    introducing a  sequence
\begin{equation}
 u_{n} (z)=  {\cal A}_{n} (z)^{-1}  f_{n} (z), \q n\in {\Bbb Z}_{+}.
\label{eq:Gs4}\end{equation}
    In proofs, we usually   omit the dependence on $z$ in notation; for example, we write $f_{n}$, $u_{n}$, ${\bf r}_{n}$.


First, we  derive a difference equation for $ u_{n} (z)$.

\begin{lemma}\label{Gs}
Let the remainder ${\bf r}_{n} (z)$ be defined by formula \e{eq:Grr}.  Set
\begin{equation}
\Lambda_{n} (z)=\frac{a_{n}}{a_{n-1} }  \frac{{\cal A}_{n+1}(z)}{{\cal A}_{n-1} (z)} 
\label{eq:GL}\end{equation}
and
\begin{equation}
{\cal R}_{n} (z) = - \sqrt{\frac{a_{n}}{a_{n-1} } }  \frac{ {\cal A}_{n}(z)}{ {\cal A}_{n-1}(z)}  {\bf r}_{n} (z).
\label{eq:GL1}\end{equation}
 Then
 equation  \e{eq:Jy} for a sequence $ f_{n} (z)$ is equivalent to the equation
\begin{equation}
\Lambda_{n} (z)( u_{n+1} (z)- u_{n} (z)) -      ( u_{n} (z)- u_{n-1} (z))=   {\cal R}_{n} (z) u_{n} (z), \q n\in {\Bbb Z}_{+},
\label{eq:DE}\end{equation}
for  sequence  \e{eq:Gs4}. 
 \end{lemma}

\begin{pf}
Substituting  expression $f_{n}  = {\cal A}_{n}  u_{n} $ into  \e{eq:Jy} and using definition \e{eq:Grr},
we see that
\begin{align*}
( \sqrt{a_{n-1}a_{n}}{\cal A}_{n} )^{-1}&\Big(  a_{n-1} f_{n-1} + ( b_{n} -z)f_{n} + a_{n} f_{n+1}\Big)
\\=&
\sqrt{\frac{a_{n-1} }{a_{n}}}\frac{{\cal A}_{n-1} }{{\cal A}_{n}}  u_{n-1} +\frac{b_{n}-z }{ \sqrt{a_{n-1}a_{n}}}  u_{n} +
\sqrt{\frac{a_{n} }{a_{n-1}}}  \frac{{\cal A}_{n+1} }{{\cal A}_{n}} u_{n+1} 
\\=&
\sqrt{\frac{a_{n-1} }{a_{n}}}\frac{{\cal A}_{n-1} }{{\cal A}_{n}} ( u_{n-1}-u_{n})+
\sqrt{\frac{a_{n} }{a_{n-1}}}  \frac{{\cal A}_{n+1} }{{\cal A}_{n}} (u_{n+1}-u_{n}) + {\bf r}_{n}  u_{n}
\\=&
\sqrt{\frac{a_{n-1} }{a_{n}}}\frac{{\cal A}_{n-1} }{{\cal A}_{n}}\Big(( u_{n-1}-u_{n})+ \Lambda_{n}( u_{n+1}  - u_{n} )  - {\cal R}_{n}u_{n}\Big)
  \end{align*}
  where the coefficients $\Lambda_{n}$ and ${\cal R}_{n}$ are defined by equalities \e{eq:GL} and \e{eq:GL1}, respectively.
 Therefore     equations \e{eq:Jy}  and \e{eq:DE} are equivalent.
 \end{pf}

   Our next goal is to construct a solution of  difference equation
\e{eq:DE} such that
  \begin{equation}
\lim_{n\to\infty} u_{n} (z)=   1.
\label{eq:DE1}\end{equation}
To that end, we will reduce    equation \e{eq:DE} to a Volterra ``integral" equation  which can be standardly solved by successive approximations.

    \subsection{Volterra equation}
    
    It is convenient to consider this problem in a more general setting. 
     We now do not make
    any specific assumptions about the sequences $\Lambda_{n} $ and $  {\cal R}_{n}  $  in \e{eq:DE} except that $\Lambda_{n}\neq 0 $.
     Denote
   \begin{equation}
  X_{n}= \Lambda_{1} \Lambda_{2}\cdots \Lambda_{n}  
\label{eq:W4}\end{equation}
and
\begin{equation}
G_{n,m}  =  X_{m-1}   \sum_{p=n }^{m-1}   X_{p}  ^{-1} , \q   m \geq n+1.
\label{eq:DE3}\end{equation}
The sequence $u_{n}$ will be constructed as a solution of
 a discrete  Volterra integral equation 
 \begin{equation}
u_{n}  = 1+ \sum_{m=n+1}^\infty G_{n,m}   {\cal R}_m  u_{m} .
\label{eq:DE6}\end{equation}

Under natural assumptions, this equation can be standardly solved by successive approximations. First, we estimate its iterations.

  \begin{lemma}\label{GS3p}
   Let us set
\begin{equation}
 h_{m} =\sup _{n\leq m-1} |G_{n,m}   {\cal R}_{m}  |
\label{eq:DE4}\end{equation}
 and suppose that
\begin{equation}
  ( h_{m} )\in \ell^1 ({\Bbb Z}_{+}).
\label{eq:DE5}\end{equation}
Put $u^{(0)}_n =1$   and 
  \begin{equation}
 u^{(k+1)}_{n}  = \sum_{m=n +1}^\infty G_{n,m}    {\cal R}_{m}  u^{(k )}_m  ,\q k\geq 0,
\label{eq:W5}\end{equation}
for all $n\in {\Bbb Z}_{+}$. Then  estimates
 \begin{equation}
| u^{(k )}_{n}   |\leq \frac{H_{n} ^k}{k!} ,\q \forall k\in{\Bbb Z}_{+},
\label{eq:W6s}\end{equation}
where
\begin{equation}
 H_n = \sum_{p=n +1}^\infty  h_p , 
\label{eq:DEs}\end{equation}
are true. 
\end{lemma}

 \begin{pf}
  Suppose that \e{eq:W6s} is satisfied for some $k\in{\Bbb Z}_{+}$. We have to check 
 the same estimate (with $k$ replaced by $k+1$ in the right-hand side)  for $ u^{(k+1)}_{n}$.  
  According to definitions \e{eq:DE4} and \e{eq:W5}, it  follows from estimate \e{eq:W6s} that
   \begin{equation}
| u^{(k +1)}_{n} |\leq  \sum_{m=n +1}^\infty  h_m | u^{(k )}_{m} |\leq \frac{1}{k!}  \sum_{m=n +1}^\infty  h_m H_{m}^k.
\label{eq:V7}\end{equation}
Since $H_{m-1}= H_{m}+ h_{m}$, we have an inequality
 \[
H_{m}^{k+1}+ (k+1)   h_{m} H_{m}^k  
\leq (H_{m}+h_{m})^{k+1}
=H_{m-1}^{k+1},
\]
and hence, for all $N\in{\Bbb Z}_{+}$,
   \[
 (k+1)  \sum_{m=n +1}^N  h_m   H_{m}^k 
 \leq 
 \sum_{m=n +1}^N  ( H_{m-1}^{k+1}- H_{m}^{k+1})
= H_{n}^{k+1}- H_{N}^{k+1}\leq H_{n}^{k+1}.
 \]
Substituting this bound into   \e{eq:V7}, we obtain estimate \e{eq:W6s} for $u^{(k +1)}_{n}$.
    \end{pf}

Now we are in a position to solve equation \e{eq:DE6} by iterations.

  \begin{theorem}\label{GS3}
Let assumption  \e{eq:DE5} be satisfied. Then equation  \e{eq:DE6}
  has a  bounded solution $u_{n} $. This solution satisfies  an estimate
  \begin{equation}
 |u _{n}   -1 |\leq e^{H_{n}}-1 \leq C H_{n}
\label{eq:Gp9}\end{equation}
where $H_{n}$ is sum  \e{eq:DEs}.
In particular, condition \e{eq:DE1} holds. 
\end{theorem}

 \begin{pf}  Set
     \begin{equation}
   u_{n} =\sum_{k=0}^\infty u^{(k)}_{n} 
\label{eq:W8}\end{equation}
where $u^{(k)}_{n}$ are defined by recurrence relations \e{eq:W5}.
Estimate \e{eq:W6s} shows that this series is absolutely convergent. Using the Fubini theorem to interchange the order of summations in $m$ and $k$, we see that
\[
   \sum_{m=n+1}^\infty G_{n,m}     {\cal R}_{m}  u_{m}    =   \sum_{k=0}^\infty\sum_{m=n+1}^\infty G_{n,m}   {\cal R}_{m}  u_{m}^{(k)} = \sum_{k=0}^\infty  u_{n}^{(k+1)}=-1+\sum_{k=0}^\infty  u_{n}^{(k)}=-1+  u_{n}.
\]
This is equation  \e{eq:DE6} for sequence \e{eq:W8}. Estimate \e{eq:Gp9}  also follows from \e{eq:W6s} and \e{eq:W8}.  
\end{pf} 

 \begin{remark}\label{GS3r}
  A bounded solution $u_{n} $ of \e{eq:DE6} is of course unique. Indeed, suppose  that $ ( v_{n} )
  \in \ell^\infty ({\Bbb Z}_{+})$ satisfies  homogeneous equation  \e{eq:DE6}, that is,
   \[
v_{n}   =  \sum_{m=n+1}^\infty G_{n,m}  {\cal R}_m   v_{m}  .
\]
Then, by assumption \e{eq:DE5}, we have
  \[
| v_{n}  |\leq  \sum_{m=n+1}^\infty h_m   |v_{m}  |.
\]
Iterating this estimate, we find that
\[
| v_{n}  |\leq    \frac{1}{k!}   H_{n}^k  \max_{n\in {\Bbb Z}_{+}} \{| v_{n}  |\}, \q \forall k \in {\Bbb Z}_{+}.
\]
Taking the limit $k\to\infty$, we see  that $v_{n}=0$. Note however that we do not use the unicity in our construction.
\end{remark}

   \subsection{Back to  the difference equation}   
 
It turns out that the construction above yields a solution of  difference equation  \e{eq:DE}.
 
  \begin{lemma}\label{GS4}
  Under assumption \e{eq:DE5} 
 a solution $u_{n}   $ of   integral  equation   \e{eq:DE6}   satisfies an identity
 \begin{equation}
   u_{n+1} -    u_{n} =  - X_{n}^{-1} \sum_{m=n+1}^\infty        X_{m-1}  {\cal R}_{m} u_{m} 
  \label{eq:A17Mb}\end{equation}
 and  difference equation  \e{eq:DE}. 
 \end{lemma}
 
  \begin{pf}
  It follows from \e{eq:DE6}    that
  \begin{equation}
 u_{n+1} - u_{n}=   \sum_{m=n+2}^\infty (G_{n+1, m}-G_{n, m}){\cal R}_m u_{m} -
 G_{n, n+1} {\cal R}_{n+1} u_{n+1}.
  \label{eq:A17Ma}\end{equation}
  According to  \e{eq:DE3}  we have
\[
G_{n+1, m}-G_{n, m}=-  X_{n}^{-1} X_{m-1} \q \mbox{and}\q
G_{n, n+1}=     1 .
\]
Therefore relation \e{eq:A17Ma} can be rewritten as
\e{eq:A17Mb}.

Putting together  equality \e{eq:A17Mb} with the same equality  where  $n+1$ is replaced by $n$, we see that
 \[
 \Lambda_{n} ( u_{n+1} - u_{n}) -   ( u_{n} - u_{n-1})=
 - \Lambda_{n} X_{n}^{-1} \sum_{m=n+1}^\infty        X_{m-1}  {\cal R}_{m} u_{m} +
 X_{n-1}^{-1} \sum_{m=n}^\infty        X_{m-1}  {\cal R}_{m} u_{m} .
  \]
Since $X_{n}=  \Lambda_{n} X_{n-1}$,   the right-hand side here equals $ {\cal R}_{n}u_{n}$, and hence the equation obtained coincides with \e{eq:DE}.  
   \end{pf}

 
    \begin{corollary}\label{DEe}
    It follows  from  \e{eq:A17Mb} that
\begin{equation}
| u_{n+1} - u_{n}|  \leq \max_{n\in {\Bbb Z}_{+}}\{ |u_{n}| \} \, |X_{n}|^{-1} \sum_{m=n}^\infty  |      X_{m}  {\cal R}_{m+1} |.
  \label{eq:A1y}\end{equation} 
   \end{corollary}

    Lemma~\ref{GS4} allows us to reformulate Theorem~\ref{GS3} in terms of solutions of equation   \e{eq:DE}.

     \begin{theorem}\label{DE}
Let assumption  \e{eq:DE5} be satisfied. Then  difference equation  \e{eq:DE}
  has a   solution $u_{n} (z)$  satisfying  estimates
\e{eq:Gp9} and \e{eq:A1y}.
In particular, condition \e{eq:DE1} holds. 
\end{theorem}

 Let us now discuss the dependence on the spectral parameter $z$.   
   Suppose   that the coefficients  $\Lambda_{n} (z)$ and ${\cal R}_{n} (z)$ in equation   \e{eq:DE6} are functions of $z\in\Omega$ on some open set $\Omega\subset{\Bbb C}$. 
   
   
    \begin{lemma}\label{DE1}
Let  the coefficients  $\Lambda_{n} (z)$ and   ${\cal R}_{n} (z)$ be analytic functions of $z\in  \Omega$. Suppose that
 assumption  \e{eq:DE5}   is  satisfied uniformly in $z $ on compact subsets of $\Omega$. Then the solutions $u_{n}  (z)$ of integral equation   \e{eq:DE6}  are also analytic in $z\in \Omega$. Moreover, if $\Lambda_{n} (z)$ and  ${\cal R}_{n} (z)$ are continuous up to the boundary of $\Omega$  and  assumption \e{eq:DE5}  is   satisfied uniformly on $\Omega$, then the same is true for the functions $u_{n}  (z)$.
\end{lemma}

  \begin{pf}
  Consider  series \e{eq:W8} for a solution   $u_{n} (z)$ of integral equation \e{eq:DE6}. 
 Observe that if the functions $u_{m}^{(k)} (z)$ in \e{eq:W5} depend analytically (continuously)  on $z$, then the function $u_{n}^{(k+1)} (z)$ is also analytic (continuous). Since  
  series  \e{eq:W8}  converges uniformly, its sums $u_{n}  (z)$ are also   analytic (continuous)  functions.  
       \end{pf}
       
        In view of Lemma~\ref{GS4} this result applies also   to solutions of difference equation  \e{eq:DE}.

     
      \section{Jost solutions}
 
 Here we use the results of the previous section to  construct the Jost solutions $f_{n}  (z)$ of the Jacobi equation \e{eq:Jy}
  with the coefficients $a_{n} $ and $b_{n}$ satisfying conditions \e{eq:ASa} and \e{eq:ASb} where $|\gamma| =1$.  This leads to Theorems~\ref{GSS-x}, \ref{GSS-y} and \ref{GSS-z}.   
  
  First, in Sect.~5.1 and 5.2, we state  some necessary technical results.

  \subsection{Discrete derivatives}


Let us collect standard formulas for  ``derivatives" 
 \[
x_{n}'= x_{n+1}  -x_{n}
\]
of  various sequences $x_{n}$:
 \begin{align}
(x_{n}^{-1})'= &- x_{n}^{-1}x_{n+1}^{-1} x_{n}',
\nonumber\\
(e^{x_{n}})'=  &(e^{ x_{n}'}-1) e^{x_{n}},
\label{eq:dife}\\
(\sqrt{x_{n }})'= &  x_{n}'  (\sqrt{x_n }+\sqrt{x_{n+1 }})^{-1}
\nonumber \end{align} 
and
 \begin{equation}
(x_{n} y_{n})'=  x_{n+1} y_{n}'+ x_{n}' y_{n}.
\label{eq:dife2}\end{equation}
 Note the Abel summation formula (``integration by parts"):
 \begin{equation}
\sum_{p=n }^m x_p y_p' = x_{m}  y_{m+1} - x_{n -1}  y_{n}  -\sum_{p=n } ^m x_{p-1}'  y_{p};
\label{eq:Abel}\end{equation}
here $m \geq n\geq 0$ are arbitrary (we   set $x_{-1}=0$ so that $x_{-1}'=x_{0}$).

We mention also an obvious estimate
  \begin{equation}
    | f (x_{n+1})-f (x_{n})|\leq \big( \max_{ |x| \leq 1}    | f' (x)| \big) | x_{n}'|
    \label{eq:kk}\end{equation}
    valid for an arbitrary function $f\in C^1$,  an arbitrary sequence $x_n\to 0$ as $n\to\infty$ and sufficiently large $n$.

  Let us now consider equation \e{eq:Jy}.  A direct calculation shows that, for two  $f=(f_{n})_{n=-1}^\infty$ and $g=(g_{n})_{n=-1}^\infty$
  solutions of this equation, their Wronskian
    \begin{equation}
   W [ f,g ]:= a_{n}  (f_{n}  g_{ n+1}-f_{ n+1}  g_{n})
\label{eq:Wr}\end{equation}
does not depend on $n=-1,0,1,\ldots$. In particular, for $n=-1$ and $n=0$, we have
   \[
   W[f,g]= a_{-1} (f_{-1}  g_{0}-f_{0}  g_{-1}) \q {\rm and} \q  W [ f,g ] = a_{0}  (f_{0}  g_{ 1}-f_{ 1}  g_{0})
\]
(the number $a_{-1}\neq 0$ is arbitrary,  but the products $a_{-1} f_{-1}$ do not depend on its choice).
Clearly, the Wronskian $W [ f,g ] =0$ if and only if the solutions $f$ and $g$ are proportional.

    \subsection{Oscillating sums}

Below we need to estimate sums of oscillating or exponentially growing terms. First, we note an integration-by-parts formula. The following elementary assertion does not require specific assumptions  about amplitudes $\boldsymbol\kappa_{n}$ and phases $\boldsymbol\varphi_{n}$.

\begin{lemma}\label{IP}
Set $\boldsymbol\theta_{n}=\boldsymbol\varphi_{n+1}-\boldsymbol\varphi_{n}$ and
  \begin{equation}
\boldsymbol\zeta_{n}= \boldsymbol\kappa_{n} (e^{-i \boldsymbol\theta_{n}}-1)^{-1}.
    \label{eq:kz}\end{equation}
Then
  \begin{equation}
    \sum_{p=n}^m  \boldsymbol\kappa_{p} e^{-i\boldsymbol\varphi_p }=   \boldsymbol\zeta_{m} e^{-i\boldsymbol\varphi_{m+1} }-
     \boldsymbol\zeta_{n-1} e^{-i\boldsymbol\varphi_n  }-  \sum_{p=n}^m \boldsymbol\zeta_{p-1}' e^{-i\boldsymbol\varphi_p }
    \label{eq:k}\end{equation}
 for all $n$ and $m$.
 \end{lemma}
 
 \begin{pf}
 According to  \e{eq:dife}  the left-hand side of \e{eq:k} can be rewritten as
 \[
  \sum_{p=n}^m  \boldsymbol\zeta_{p} ( e^{-i\boldsymbol\varphi_p })'.
  \]
   It follows from formula  \e{eq:Abel} that  this sum   equals the right-hand side of \e{eq:k}.
 \end{pf}
 
 \begin{corollary}\label{IP1}
 Suppose that 
   \begin{equation}
  \boldsymbol\zeta_{n}'  \in \ell^1 ({\Bbb Z}_{+})
    \label{eq:kx}\end{equation}
     and $\Im \boldsymbol\theta_{n}\geq 0$.
Then  
  \begin{equation}
   \big| \sum_{p=n}^m   \boldsymbol\kappa_{p} e^{-i\boldsymbol\varphi_p }\big| \leq C e^{\Im\boldsymbol\varphi_{m+1} }
    \label{eq:k1}\end{equation}
 where the constant $C$ does not depend on $n$ and $m$.
 \end{corollary}
 
  \begin{remark}\label{IP2}
  If 
    \begin{equation}
\boldsymbol\theta_{n}'  \in \ell^1 ({\Bbb Z}_{+}),
    \label{eq:kx1}\end{equation}
 then condition \e{eq:kx}  can be replaced by  more convenient ones:
  \begin{equation}
  \frac{\boldsymbol\kappa_{n}} {\boldsymbol\theta_{n}} \in \ell^\infty ({\Bbb Z}_{+})\q \mbox{and}\q
\big(\frac{\boldsymbol\kappa_{n}} {\boldsymbol\theta_{n}}\big)' \in \ell^1 ({\Bbb Z}_{+}).
   \label{eq:k2}\end{equation}
 \end{remark}

   \begin{pf}
   It follows from \e{eq:dife2} that
  \begin{equation}
   \boldsymbol\zeta_{n}' = \big(\frac{\boldsymbol\kappa_{n}} {\boldsymbol\theta_{n}}\big)'
  \frac{\boldsymbol\theta_{n+1}} {e^{-i\boldsymbol\theta_{n+1}}-1} +  \frac{\boldsymbol\kappa_{n}} {\boldsymbol\theta_{n}} 
  \big( \frac{\boldsymbol\theta_{n}} {e^{-i\boldsymbol\theta_{n}}-1}\big)'.
  \label{eq:k3}\end{equation}
  Note that the function $f (\boldsymbol\theta)= \boldsymbol\theta (e^{-i\boldsymbol\theta}-1)^{-1}$ is $C^1$ in a neighbourhood of the point $\boldsymbol\theta=0$. Therefore the sequence $f (\boldsymbol\theta_{n})$ is bounded as $n\to\infty$ and  $f' (\boldsymbol\theta_{n})\in\ell^1 ({\Bbb Z}_{+})$ according to  estimate \e{eq:kk} and condition   \e{eq:kx1}. Thus  conditions \e{eq:k2}  imply that both terms in the right-hand side of \e{eq:k3} are in $\ell^1 ({\Bbb Z}_{+})$. 
      \end{pf}

        \subsection{Estimate of the ``integral" kernel}
     
  Recall that the sequences ${\cal A}_{n}={\cal A}_{n}(z)$ and $\Lambda_{n} = \Lambda_{n}  (z)$ are  given  by relations 
     \e{eq:Gr1} and \e{eq:GL}, respectively.    Our goal  is to estimate  the matrix elements $G_{n,m}  $  defined by equalities \e{eq:W4} and \e{eq:DE3} 
 and to prove inclusion \e{eq:DE5}. 
   Our estimates apply to all values of $\sigma$.

     Putting together formulas \e{eq:GL} and \e{eq:W4}, we see that
       \[
X_{n}= c a_{n}{\cal A}_{n+1} {\cal A}_{n}
\]
where the constant $c=(a_{0} {\cal A}_{1} {\cal A}_{0})^{-1}$.
       According to  definition   \e{eq:Gr1} this yields  equality 
      \begin{equation} 
 X_{n}^{-1}=  -  c \boldsymbol\kappa_{n}  e^{- i \boldsymbol\varphi_{n}}
\label{eq:nu1}\end{equation}
where
  \begin{equation} 
  \boldsymbol\kappa_{n}= n^{\rho} (n+1)^{\rho}a_{n}^{-1}
\label{eq:nu2}\end{equation}
and
  \begin{equation} 
 \boldsymbol\varphi_{n}=\varphi_{n} + \varphi_{n+1}.
\label{eq:nu3}\end{equation}
It follows from condition  \e{eq:ASa}  that
  \begin{equation} 
 {\boldsymbol\kappa}_{n} = n^{-\nu} \big(1+ (\rho-\alpha) n^{-1} + O (n^{-2}) \big)
  \label{eq:nu3A}\end{equation}
where $\nu=\sigma - 2\rho$ satisfies \e{eq:nu}.



     First, we reformulate Lemma~\ref{IP}  and its consequences in a particular form adapted to our problem.
     
     
\begin{lemma}\label{sum}
Let the assumptions of one of Theorems~\ref{GSS-x}, \ref{GSS-y} or \ref{GSS-z} be satisfied.  Define
 the sequences  $  \boldsymbol\kappa_{n}$ and  $  \boldsymbol\varphi_{n}$  by equalities \e{eq:nu2}  and \e{eq:nu3}. 
Then  estimate \e{eq:k1} holds.
 \end{lemma}

 \begin{pf}
 Set  
     \begin{equation}
   \boldsymbol\theta_{n}=   \boldsymbol\varphi_{n+1}-  \boldsymbol\varphi_{n}=\theta_{n}+\theta_{n+1}.
      \label{eq:kv}\end{equation}
 It follows from relations \e{eq:hold} or \e{eq:hold11}  that inclusion \e{eq:kx1} holds.  Therefore in view of Remark~\ref{IP2},   it suffices to check inclusions \e{eq:k2}.  By definition  \e{eq:bt1-}, we have
    \begin{equation}
 \boldsymbol\kappa_{n}\boldsymbol\theta_{n}^{-1}=  \big(n^{\nu}\boldsymbol\kappa_{n} \big)\big(n^{\nu} S_{n}\big)^{-1}, \q S_{n}=\sqrt{T_{n}}+ \sqrt{T_{n+1}},
   \label{eq:k2x}\end{equation}
   where $T_{n}$ is defined by equality \e{eq:te} (in particular, $T_{n}=t_{n}$ if $\sigma>2/3$).  
    Inclusions $  \big(n^{\nu}\boldsymbol\kappa_{n} \big) \in {\ell}^\infty ({\Bbb Z}_{+})$  and $  \big(n^{\nu}\boldsymbol\kappa_{n} \big)' \in {\ell}^1({\Bbb Z}_{+})$ are direct
   consequences of formula \e{eq:nu3A}. It follows from  relations \e{eq:bt2},  \e{eq:bt2+}  or  \e{eq:hm}
   that the product 
  $ n^{\nu}  S_{n}$ has a finite non-zero limit (it is used here that  $z\neq 0$ under the assumptions of   Theorem  \ref{GSS-y}  and that $z\neq \gamma \tau $ under the assumptions of   Theorem  \ref{GSS-z}).  The inclusion $ \big(n^{\nu} S_{n}\big)' 
   \in {\ell}^1({\Bbb Z}_{+})$ is again a consequence of \e{eq:hold} or \e{eq:hold11}.
  Therefore \e{eq:k2x} implies inclusions \e{eq:k2} which yields \e{eq:k1}. 
   \end{pf}

 Now we are in a position to estimate the matrix elements $G_{n,m}  $.  First, we note that
    \begin{equation}
 C_{1} m^{\nu} e^{ - \Im  \boldsymbol\varphi_m }  \leq |X_{m}| \leq C_{2} m^{\nu} e^{ - \Im  \boldsymbol\varphi_m } 
\label{eq:nu5}\end{equation}
according to definition  \e{eq:nu1}  and  relation  \e{eq:nu3A}.
 Next, we apply inequality  \e{eq:k1} to   elements \e{eq:nu1} which yields
   \begin{equation}
\big| \sum_{p=n }^{m-1}   X_{p}  ^{-1} | \leq  C e^{\Im  \boldsymbol\varphi_m } .
\label{eq:nu41}\end{equation}
Combining \e{eq:nu5} and \e{eq:nu41}, we obtain a convenient 
  estimate on product \e{eq:DE3}.

\begin{lemma}\label{GS2+}
Under the assumptions of any of Theorems~\ref{GSS-x}, \ref{GSS-y} or \ref{GSS-z}, we have an estimate
\begin{equation}
|G_{n,m}|\leq    Cm^{\nu}
\label{eq:AbG}\end{equation}
where $\nu$ is given by  \e{eq:nu} and 
the constant  $C$  does not depend on $n$ and $m$.
 \end{lemma}
 
 
\subsection{Solutions of the integral equation}

Next, we consider  integral equation \e{eq:DE6}. Observe that   remainder \e{eq:GL1}   obeys the same estimate
 \e{eq:rem} as ${\bf r}_{n}$. Thus, according to the results of Sect.~3  (see Propositions~\ref{AA+}, \ref{AA+1}, \ref{AA+-}  and \ref{AA-})
   \begin{equation}
|R_{m}|\leq    Cm^{-\d}
\label{eq:Rr}\end{equation} 
where $\d$  satisfies conditions \e{eq:R1}.

Putting together  \e{eq:AbG} and \e{eq:Rr}, we obtain an estimate on sequence \e{eq:DE4}:
  \[
h_{m} \leq    Cm^{\nu -\d}. 
\]
Comparing \e{eq:nu} and \e{eq:R1}, we see that $ \nu -\d  < -1$.  It follows that sum  \e{eq:DEs} satisfies  an estimate
  \[
H_{n}\leq    Cn^{\nu -\d +1}. 
\]
Therefore condition  \e{eq:DE5} holds, and Theorem~\ref{DE} applies in our case. This yields estimates \e{eq:Gp9} and \e{eq:A1y}. Moreover, the right-hand side of   \e{eq:A1y} can be estimated explicitly. Indeed, note that $\Im \boldsymbol\varphi_n \leq \Im \boldsymbol\varphi_m$ for $m\geq n$ according to \e{eq:the}. Thus, it follows from  \e{eq:nu1} and \e{eq:nu3A} that
\[
| X_{n}^{-1 }X_{m}| \leq  | k_{n} k_{m}^{-1 }| e^{ \Im( \boldsymbol\varphi_n -  \boldsymbol\varphi_m) } \leq C n^{-\nu } m^{\nu  } ,\q m\geq n, 
\]
so that inequality   \e{eq:A1y} yields an estimate
 \[
| u_{n+1}- u_{n}| \leq C n^{-\nu} \sum_{m= n}^\infty  m^{\nu- \d}\leq C_{1} n^{ - \d+1}.
 \]

We see that, under the assumptions of any of Theorems~\ref{GSS-x}, \ref{GSS-y} or \ref{GSS-z},   condition \e{eq:DE5} is satisfied.  Hence
the following three results are direct consequences  of Theorem~\ref{DE} (see also Lemma~\ref{DE1}). 
Recall that the number $\nu$  is defined by relations \e{eq:nu}  and   $\d$ satisfies conditions \e{eq:R1}.



 
 \begin{theorem}\label{HSS-x} 
       Let the assumptions  of Theorem~\ref{GSS-x}  be satisfied.

If $\tau<0$, then
    for every  $z\in \clos\Pi$    equation \e{eq:DE6} has a solution $  u_{n}( z ) $ satisfying asymptotic relations 
    \begin{equation}
 u_{n}( z ) =1+ O (n^{\nu - \d +1})
\label{eq:DX}\end{equation}
and
  \begin{equation}
 u_{n}' ( z ) = O (n^{- \d +1})
\label{eq:DXD}\end{equation}
where $\nu=1/2$ and $\d>1/2+\sigma$.
For all $n\in {\Bbb Z}_{+}$, the functions $u_{n}( z )$ are analytic in $\Pi$ and are continuous up to the cut along the real axis.  

If $\tau > 0$, then relations  \e{eq:DX} and  \e{eq:DXD} are true for all $z\in {\Bbb C}$. In this case  the functions $u_{n}( z )$ are analytic in the whole complex plane ${\Bbb C} $.

For all $\tau\neq 0$, asymptotic formula \e{eq:DE1} is uniform in $z$ from compact subsets    of $  \Bbb C $.
\end{theorem}


   \begin{theorem}\label{HSS-y}
           Let the assumptions  of Theorem~\ref{GSS-y}  be satisfied.  Then for  every $z\neq 0$ such that $z\in \gamma \clos \Pi_{0}$,     equation \e{eq:DE6} has a solution $  u_{n}( z ) $ with asymptotics
\e{eq:DX}, \e{eq:DXD} where $\nu=\sigma/2$ and $\d>1+\sigma/2$. 
For all $n\in {\Bbb Z}_{+}$, the functions $u_{n}( z )$ are analytic in $z\in \gamma\Pi_{0} $ and are continuous up to the cut along the   half-axis $\gamma{\Bbb R}_{+}$, with a possible exception of the point $z=0$.  
         \end{theorem}

 \begin{theorem}\label{HSS-z}
             Let the assumptions  of Theorem~\ref{GSS-z}  be satisfied. 
     Then for  every $z$ such that $ z\in \gamma (\tau+\clos \Pi_{0})$, $z\neq \gamma\tau$,   equation \e{eq:DE6} has a solution $ u_{n}( z )$ with asymptotics  \e{eq:DX},  \e{eq:DXD} where $\nu=1/2$ and $\d= 2$.
For all $n\in {\Bbb Z}_{+}$, the functions $u_{n}( z )$ are analytic in $z\in \gamma (\tau+ \Pi_{0})$  and are continuous up to the cut along the   half-axis $\gamma(\tau+ {\Bbb R}_{+})$, with a possible exception of the point $\gamma\tau$.  
              \end{theorem}


  \subsection{The Jost solutions}
                 
                   Now it is easy to construct solutions of  the Jacobi equation \e{eq:Jy} with asymptotics \e{eq:A22G+} as $n\to\infty$.
  We call them the Jost solutions. 
  
According to Lemma~\ref{Gs}  equation \e{eq:DE} for the sequence $u_{n}(z)$ and  equation \e{eq:Jy} for the sequence  
       \begin{equation}
f_{n}(  z )  =  (-\gamma)^n n^{- \rho} e^{  i \varphi_{n} (\gamma z)} u_{n} (z)
\label{eq:wrua}\end{equation}
are equivalent.  Therefore Theorems~\ref{GSS-x}, \ref{GSS-y} and \ref{GSS-z} are direct consequences of Theorems~\ref{HSS-x}, \ref{HSS-y} and \ref{HSS-z},   respectively.


 
 
 Finally, we show that the Jost solutions are determined uniquely by their asymptotics \e{eq:A22G+}. 
   This is quite simple for regular $z$.  Recall that the set ${\cal S}$  was defined by relations \e{eq:SS}. 
 
  \begin{proposition}\label{uniq}
             Let the assumptions  of one of Theorems~\ref{GSS-x}, \ref{GSS-y}  or \ref{GSS-z}   be satisfied.   If  $\sigma=3/2$, we also assume that $\tau>0$.
      Suppose  that $z\not\in \clos{\cal S}$.       Then the solution of $f_{n} (z)$ of equation 
            \e{eq:Jy}       satisfying condition    \e{eq:A22G+}    is unique. 
                \end{proposition}
                
                 \begin{pf}
                 Suppose that solutions $f_{n}$ and $\ti{f}_{n}$ of equation      \e{eq:Jy} are given by equality \e{eq:wrua} where $u_{n}$ and $\ti{u}_{n}$ obey condition        \e{eq:DE1}. Then their Wronskian  \e{eq:Wr} equals
                           \begin{equation}
W[ f,\ti{f}]=-\gamma a_{n}n^{-\rho} (n+1)^{-\rho}e^{i\varphi_{n}}e^{i\varphi_{n+1}}\big(u_{n}\ti{u}_{n+1}-u_{n+1}\ti{u}_{n}\big).
\label{eq:wru}\end{equation}
As explained in Sect.~2.2, under  the assumptions of Proposition~\ref{uniq} the sequence $e^{i\varphi_{n} }$  tends to zero exponentially as $n\to\infty$ whence 
 $W[ f,\ti{f}]= 0 $ and consequently
 $\ti{f}_{n}=C f_{n}$ for some constant $C$.   It now  follows  from \e{eq:wrua} that    $\ti{u}_{n}=C u_{n}$  where  $C=1$ by  \e{eq:DE1}.
                                 \end{pf}
                                 
                                   \begin{remark}\label{uniq-r}
                                    If  $\sigma=3/2$ and $\tau<0$, then  instead of \e{eq:DE1} we have to require a stronger condition 
                                           \begin{equation}
                                           u_{n} =1+O(n^{-1/2}). 
                                           \label{eq:wruA}\end{equation}
                                           Note that in view of \e{eq:Af1-} this condition  is satisfied for the Jost solution $f_{n} (\lambda+i\varepsilon)$  of equation      \e{eq:Jy} constructed in Theorem~\ref{GSS-x}.  Suppose that two solutions $f_{n}$ and $\ti{f}_{n}$ are given by formula   \e{eq:wrua} where $u_{n}$ and $\ti{u}_{n}$ satisfy \e{eq:wruA} whence $u_{n}\ti{u}_{n+1}-u_{n+1}\ti{u}_{n}= O(n^{-1/2})$. Since $\rho=1/2$ now, it follows from asymptotic formula  \e{eq:Af1-}  and relation                                    \e{eq:wru} that
                                         \[
|W[ f,\ti{f}]|= O \big(a_{n}n^{-1- 2|\varepsilon|/  \sqrt{|\tau|} } (u_{n}\ti{u}_{n+1}-u_{n+1}\ti{u}_{n} )\big) = O (n^{ -2|\varepsilon|/  \sqrt{|\tau|} })
= 0 
\]
because $\varepsilon\neq 0$.  This implies that $\ti{f}_{n}=f_{n}$. 
                \end{remark}

                                 The results for $z$ in the spectrum of the operator $J$ are slightly weaker.
                                 
                                   \begin{proposition}\label{uniq+}
             Let the assumptions  of one of Theorems~\ref{GSS-x}, \ref{GSS-y}  or \ref{GSS-z}   be satisfied. Suppose  that  $z=\lambda\pm i0$
             where $\lambda\in {\cal S}$.
             Then the solution $f_{n} (z)$ of equation 
            \e{eq:Jy}       satisfying relation     \e{eq:wrua}  with $u_{n}$ obeying conditions \e{eq:DE1}  and  \e{eq:DXD} is unique. 
                \end{proposition}
                
                  \begin{pf}
                 Suppose that two solutions $f_{n}$ and $\ti{f}_{n}$ of equation      \e{eq:Jy} satisfy these conditions. Their Wronskian is given by equality                  \e{eq:wru}                where
                 \[
                 u_{n}\ti{u}_{n+1}-u_{n+1}\ti{u}_{n} =          u_{n}(  \ti{u}_{n +1}-  \ti{u}_{n}) + (u_{n}-u_{n+1})\ti{u}_{n}=O(n^{-\d +1}).
                 \]
It follows that
                    \[
W[ f,\ti{f}] =  O(n^{\nu-\d +1}), \q \nu= \sigma-2\rho, \q n\to\infty.
\]
Putting together relations  \e{eq:nu} and  \e{eq:R1}, we see  that
$
\nu-\d +1 < 0
$
for all $\sigma\in (0,3/2]$.
Therefore $W[ f,\ti{f}] =0$ and, consequently,  $\ti{f} =f$.
   \end{pf}

  \section{Orthogonal polynomials}
    
      Here we   describe an asymptotic behavior as $n\to\infty$ of all solutions $F_{n}(z)$ of equation \e{eq:Jy}. In particular,  these results apply to the orthonormal polynomials $P_{n}(z)$.  We have to distinguish values of $z=\lambda \in{\cal S}$ (this set was defined by relations   \e{eq:SS}) in the absolutely continuous  spectrum of a Jacobi operator and  regular points $z \in {\Bbb C}\setminus \clos {\cal S} $.

 \subsection{Regular points}
 
      
         
         Our goal in this subsection is to prove Theorem~\ref{P2+} and Proposition~\ref{D1+}.
 Let us proceed from the following assertion.
  
       \begin{proposition}\label{D}\cite[Theorem~2.2]{inf}
       Let $f(z)=(f_{n} (z))$ be an arbitrary solution of the Jacobi equation \e{eq:Jy} such that $f_{n} (z)\neq 0$ for sufficiently large $n$, say $n\geq n_{0}$.  Then sequence  $g(z)=(g_{n} (z))$ defined by \e{eq:GEg} also satisfies  equation \e{eq:Jy}, and the Wronskian
       \[
       W[f(z), g(z)]=1,
       \]
       so that the solutions $f(z)$ and $g(z)$ are linearly independent.
 \end{proposition}
 
 In this subsection,  we suppose that  $z \in {\Bbb C}\setminus \clos {\cal S} $ and $f_{n}= f_{n}(z)$ is the Jost solution of equation \e{eq:Jy}.  Its asymptotics is given by formulas \e{eq:Gr1}, \e{eq:Gr1x}. 
 Our aim is to find an asymptotic  behavior of  the solution $g_{n}  = g_{n} (z)$ as $n\to\infty$. The dependence on $z$ will be omitted in notation.  Let us  set
      \begin{equation} 
\Sigma_{n}   =  \sum_{m=n_{0}}^n (a_{m-1} f_{m-1}  f_{m} )^{-1},\q n\geq n_{0} ;
\label{eq:GEf}\end{equation}
   then   \e{eq:GEg}  reads as
         \begin{equation} 
      g_{n}  = f_{n}  \Sigma_{n}  .
\label{eq:gg}\end{equation}

Using equalities \e{eq:Gr1}, \e{eq:Jost} and notation  \e{eq:nu2}, \e{eq:nu3}, we can rewrite  sum  \e{eq:GEf} as
     \begin{equation} 
\Sigma_{n} =-\gamma  \sum_{m=n_{0}-1}^{n-1}
  \boldsymbol\kappa_{m}  {\bf u}_{m}  e^{-i  \boldsymbol\varphi_{m}  }\q \mbox{where} \q {\bf u}_{m}  = (u_{m}  u_{m+1} )^{-1}.
\label{eq:GEf1}\end{equation}
In view of identity \e{eq:dife}, we have
\[
e^{-i   \boldsymbol\varphi_{m} }=\big(e^{-i   \boldsymbol\theta_{m}}-1\big)^{-1} \big(e^{-i   \boldsymbol\varphi_{m}}\big)' ,
\]
with $\boldsymbol \theta_{m}$ given by \e{eq:kv}.
This allows us to integrate by parts  in \e{eq:GEf1}. Indeed, using formula \e{eq:Abel}, we find that 
  \begin{equation} 
-\gamma \Sigma_{n} e^{i  \boldsymbol\varphi_{n} }    =   \boldsymbol\zeta_{n-1}  {\bf u}_{n-1}   
 -  \boldsymbol\zeta_{n_{0}-2}  {\bf u}_{n_{0}-2}    e^{-i  \boldsymbol\varphi_{n_{0}-1} } e^{i  \boldsymbol\varphi_{n} }  
+   \wt{\Sigma}_{n}  e^{ i  \boldsymbol\varphi_{n} }  
\label{eq:GEf2}\end{equation}
where $ \boldsymbol\zeta_{n} $  is defined by equality  \e{eq:kz} and
   \begin{equation} 
\wt{\Sigma}_{n}=-  \sum_{m=n_{0}-1}^{n-1}
\big(  \boldsymbol\zeta_{m-1}    {\bf u}_{m-1}   \big)' e^{-i  \boldsymbol\varphi_{m} }.
\label{eq:Gg}\end{equation} 
  
 We will see that asymptotics of $\Sigma_{n}$ as $n \to\infty$ is determined by the first term in the right-hand side of expression \e{eq:GEf2}. Let us calculate it.
 Recall that   ${\bf u}_{n} \to 1$ as $n\to\infty$ according to Theorem~\ref{GS3}. Therefore
 putting together asymptotic formulas  \e{eq:wra} for
 $ \boldsymbol \theta_{n} $ and \e{eq:nu3A} for
 $ \boldsymbol\kappa_{n} $, we find  that
   \begin{equation} 
   \lim_{n\to\infty}  \boldsymbol\zeta_{n}  {\bf u}_{n}    =
\lim_{n\to\infty}  \boldsymbol\zeta_{n}     = i \varkappa
\label{eq:GEf3}\end{equation}
with the coefficient $ \varkappa =\varkappa (z)$ given by \e{eq:wr}.

The  second term in the right-hand side of  \e{eq:GEf2}  tends to zero as $n\to\infty$ due to the factor $ e^{ i  \boldsymbol\varphi_{n} }$.
The same is true for the third term.  To show this, we need to estimate  the derivatives  in  \e{eq:Gg}.  

 \begin{lemma}\label{DSa} 
 Let the sequence $\boldsymbol\zeta_{n}$ be defined by equality \e{eq:kz}. Then  
   \begin{equation}
 \boldsymbol\zeta_{n}'=O(n^{-1-\varepsilon})  
 \label{eq:zeta}\end{equation}
 for some $\varepsilon>0$.
       \end{lemma}

       \begin{pf}
  Let us write $\boldsymbol\zeta_{n}$ as  a product
       \begin{equation} 
\boldsymbol\zeta_{n} =  ( \boldsymbol\kappa_{n}  n^{\nu}) ( n^{\nu} \boldsymbol\theta_{n})^{-1} ( \boldsymbol\theta_{n} \big(e^{-i   \boldsymbol\theta_{n}}-1\big)^{-1})  , \q \nu=\sigma-2\rho,  
\label{eq:GEf3L}\end{equation}
and
    estimate all factors  separately. It follows from relation \e{eq:nu3A}  that the product $\boldsymbol\kappa_{n}  n^{\nu}$ tends to $1$ and its derivative is $O  (n^{-2})$ as $n\to\infty$.  Next, we consider $( n^{\nu} \boldsymbol\theta_{n})^{-1} $. According to definitions \e{eq:te}  and \e{eq:bt1-} we have
   \begin{equation} 
 n^{\nu} \theta_{n} = ( n^{\nu} \sqrt{t_{n}} )\sqrt{1+ \sum_{l=1}^{L-1} p_{l+1} t_{n}^l}.
\label{eq:kt}\end{equation}
By definition  \e{eq:ME1}, the factor $  n^{\nu} \sqrt{t_{n}} $ has a finite non-zero limit as $n\to\infty$. Moreover, its derivative is $O (n^{-\sigma})$ for $\sigma>1$ and  $O (n^{\sigma-2})$ for $\sigma<1$ (it is zero if $\sigma=1$).  Similarly, the derivative of the second factor in \e{eq:kt} is $O (n^{-2})$ for $\sigma>1$ and  $O (n^{-1-\sigma})$ for $\sigma<1$.  These arguments also show that 
$\theta_{n}'=O (n^{-3/2})$ for $\sigma>1$ and  $\theta_{n}'=O (n^{-\sigma/2-1})$  for $\sigma<1$.  Therefore the derivative of the third factor in \e{eq:GEf3L} is also $O (n^{-1-\varepsilon})$,  $\varepsilon>0$.  This proves estimate \e{eq:zeta} on product  \e{eq:GEf3L}. 
\end{pf}


  To estimate sum \e{eq:Gg}, we use   the following elementary assertion of a general nature.
  
   \begin{lemma}\label{MM} 
   Suppose that a sequence $x_{n}\in \ell^1 ({\Bbb Z}_{+})$ and  a sequence  $\vartheta_{n}\geq 0$.
   Set  
   \begin{equation} 
\phi_{n}=\sum_{m=0}^n \vartheta_{m}
\label{eq:MM}\end{equation}
and assume that 
 \begin{equation} 
  \lim_{n\to\infty}  \phi_{n} =\infty.
\label{eq:MM1}\end{equation}
Then  
         \[
     \lim_{n\to\infty}       e^{-\phi_{n}}
         \sum_{m=0}^n     x_{m} e^{\phi_{m}} =0.
\]
      \end{lemma}
      
                \begin{pf}
                By definition 
          \e{eq:MM}, we have
                      \begin{equation} 
         e^{-\phi_{n}}
         \sum_{m=0}^n     x_{m} e^{\phi_{m}} = 
         \sum_{m=0}^\infty     X_{m} (n)  
\label{eq:MM3}\end{equation}
where
              \[
           X_{m} (n) =  x_{m} \exp\big(-\sum_{p=m}^n \vartheta_{p} \big)\q \mbox{if} \q m\leq n
\]
and
$X_{m} (n)=0$ if $m>n$. Clearly, $X_{m} (n)\leq x_{m}$ because   $\vartheta_{n}\geq 0$ and $X_{m} (n)\to 0$ as $n\to\infty$  for fixed $m$ by virtue of condition \e{eq:MM1}. Therefore, by the dominated convergence theorem,  sum  \e{eq:MM3}  tends to zero as $n\to\infty$.
          \end{pf}

Now we are in a position to estimate the third term in  \e{eq:GEf2}. 

 \begin{lemma}\label{DSff} 
 Sum   \e{eq:Gg} satisfies the condition
         \begin{equation} 
\lim_{n\to\infty} \wt{\Sigma}_{n} e^{i  \boldsymbol\varphi_{n} }  =   0.
\label{eq:GEs}\end{equation}
  \end{lemma}

 \begin{pf}
 It follows from estimates  \e{eq:DXD} and \e{eq:zeta} that 
      \begin{equation} 
 |( \boldsymbol\zeta_{n}    {\bf u}_{n}   )' | \leq  | \boldsymbol\zeta_{n}' |  |   {\bf u}_{n}   | +  |\boldsymbol\zeta_{n+1} |  |   {\bf u}_{n} '  | \leq C n^{- \d +1}
 \label{eq:NN}\end{equation}
 where the value of $\d$ is indicated in Theorems~\ref{HSS-x}, \ref{HSS-y} and \ref{HSS-z}.
 Therefore, by definition \e{eq:Gg} and the differentiation  formula  \e{eq:dife}, we have
   \begin{equation} 
|\wt{\Sigma}_{n} | \leq C \sum_{m=n_{0}-1}^{n-1}  m^{-\d +1}  e^{ \phi_{m}} = C  \sum_{m=n_{0}-1}^{n-1}   y_{m}\big(e^{ \phi_{m}}\big)',\q\phi_{m}= \Im \boldsymbol\varphi_m,
 \label{eq:NN1}\end{equation}
 where
  \begin{equation} 
 y_{m} =  m^{-\d+1}  
\big(e^{ \vartheta_m}-1\big)^{-1}, \q \vartheta_{m}= \phi_{m+1} -\phi_{m}.
 \label{eq:NN2}\end{equation}
  Using  relation  \e{eq:Abel} and integrating in the right-hand side of  \e{eq:NN1} by parts, we find that 
   \begin{equation} 
|\wt{\Sigma}_{n} | \leq C\Big(  y_{n-1}e^{  \phi_n}
-y_{n_{0}-2}e^{ \phi_{n_{0}-1} }
- \sum_{m=n_{0}-1}^{n-1} y_{m-1}'
    e^{\phi_m} \Big) .
 \label{eq:NN3}\end{equation}
 
 Let us  estimate expression \e{eq:NN2}.  It follows from relations  \e{eq:bt2},  \e{eq:bt2+} and  \e{eq:hm}  that
   \[
 \phi_n= c n^{-\mu}(1+ o(1))
 \]
 for some $c= c_{\sigma,\tau}>0$. Here $\mu= \sigma/2$ if $\sigma\leq 1$,  $\mu= 1/2$ if $\sigma\in [1, 3/2]$, $\tau>0$  and  $\mu = \sigma- 1/2$ if $ \sigma\in [1, 3/2]$, $\tau<0$.  Therefore product \e{eq:NN2} is estimated as
   \begin{equation} 
| y_{n} | \leq C  n^{-\d +1+\mu}   .
 \label{eq:NN5}\end{equation}
 Note that  $-\d +1+\mu < 0$ for all values of $\sigma$ and $\tau$. Moreover,  estimate \e{eq:NN5} can be differentiated which yields
    \[
| y_{n} | \leq C  n^{ -\varepsilon} , \q  
| y_{n}' | \leq C  n^{-1 -\varepsilon} 
\]
 for some $\varepsilon>0$.

 Thus, it follows from  inequality \e{eq:NN3}  that
    \[
e^{ - \phi_n}|\wt{\Sigma}_{n} | \leq C\Big( n^{-\varepsilon}
 +\sum_{m=n_{0}-1}^{n-1} m^{-1-\varepsilon}
    e^{-\phi_n+\phi_m} \Big) 
 \]
 which in view of  Lemma~\ref{MM}  implies relation  \e{eq:GEs}.  
   \end{pf}
 
   Let us now recall  equality \e{eq:GEf2} and put relations \e{eq:GEf3} and \e{eq:GEs}  together. This leads to the following result.
   
          \begin{lemma}\label{D1} 
           Sum   \e{eq:GEf} satisfies the condition
            \begin{equation} 
\lim_{n\to\infty}  \Sigma_{n} e^{i  \boldsymbol\varphi_{n} }  = -i\gamma  \varkappa.
\label{eq:GEf4}\end{equation}
      \end{lemma}

    Using equality \e{eq:gg}   we can now conclude the {\it proofs  } of Theorem~\ref{P2+} and Proposition~\ref{D1+}. Indeed, combining asymptotics \e{eq:A22G+} and \e{eq:GEf4}, we obtain relation     \e{eq:A2P3}. This implies both formulas \e{eq:Af1-} and \e{eq:gg+}. $\q \Box$     
    
     Recall (see Sect.~1.1, for more details)  that equation  \e{eq:Jy}  is in the limit point case if, for $\Im z\neq 0$, it has a unique, up to a constant factor, non-trivial solution $f_{n} (z)$  such that 
 inclusion \e{eq:Af1+} is satisfied.  This is equivalent to the essential self-adjointness of  the minimal Jacobi operator $J_{\rm min}$  in the space $\ell^2 ({\Bbb Z}_{+})$. In this case we set $\clos J_{\rm min}=J_{\rm max}=:J$.

    According to Theorem~\ref{P2+} for $\Im z\neq 0$, the sequences $g_{n} (z)$ tend to infinity exponentially as $n\to\infty$ and according to Proposition~\ref{D1+}  they tend to infinity  as a power of $n$   (or to zero but slower than $n^{-1/2}$). In all cases relation \e{eq:gL_2} is satisfied. Therefore it follows from the limit point/circle  theory that under our assumptions the  operators $J_{\rm min}$ are essentially self-adjoint. This proves Proposition~\ref{Self-Adj1}.

      Now it is easy find an asymptotics of all solutions $F= (F_{n})$ of equation \e{eq:Jy}. Indeed, using Proposition~\ref{D}, we see that
      \[
      F_{n} = -   W[ F, f] g_{n} + c f_{n}
      \]
      for some constant $c$. The asymptotics of the solutions $g_{n} $ and $ f_{n}$ are given by formulas \e{eq:A2P3} and \e{eq:A22G+}. Obviously,  $ f_{n}$  makes no contribution to the asymptotics of $ F_{n}$. This leads to the following result.

     \begin{theorem}\label{P2}
        Let one of the following three assumptions be satisfied:
        
       $1^0$ the conditions of Theorem~\ref{GSS-x} where either $\tau<0$ and $\Im z \neq 0$ or 
        $\tau>0$ and  $z\in{\Bbb C}$ is arbitrary
       
       $2^0$ the  conditions of  Theorem~\ref{GSS-y} where either $\gamma>0$ and $ z\not\in [0,\infty)$ or 
        $\gamma<0$ and $ z\not\in (-\infty, 0]$

       $3^0$   the conditions of  Theorem~\ref{GSS-z} where either $\gamma>0$ and $ z\not\in [\tau,\infty)$ or 
        $\gamma<0$ and $ z\not\in (-\infty, -\tau]$

Then   an arbitrary   solution $ F(z) =(F_{n} (z)) $  
  has an asymptotics, as $n\to \infty$, 
            \[
F_{n}(  z )  = -  i W[ F(z), f(z)] \varkappa (z)(-\gamma)^{n+1} n^{- \rho} e^{-i\varphi_{n} (\gamma z)}\big(1 + o( 1)\big) , \q z\not\in\clos{\cal S},
\]
where the coefficient $ \varkappa   (z)$  is given by formula  \e{eq:wr}.
 \end{theorem}

 In particular, Theorem~\ref{P2}  applies to the orthonormal polynomials  $P_{n}(  z )$.  Apparently, in the critical case $|\gamma|=1$,  an asymptotic behavior of  the orthonormal  polynomials $P_{n}(z)$  for regular  points $z\in{\Bbb C}$    was never investigated before (except of  the Laguerre polynomials).  This is technically the most difficult part of this paper.

  
     \subsection{Continuous  spectrum}
      
       

  First, we check that, on the continuous spectrum of the operator $J$,  the Jost solutions $f_{n}(\lambda+i0)$ and $f_{n}(\lambda-i0) = \ov{f_{n}(\lambda+i0)}$ of equation   \e{eq:Jy} are linearly independent.  Recall that the Wronskian of two solutions of this equation  is given by formula  \e{eq:Wr}, 
        the number $\rho $ is  defined by equalities \e{eq:rho} and  the sequences $\theta_{n} (\lambda)$, $\varphi_{n} (\lambda)$ are constructed in Theorems~\ref{GSS-x}, \ref{GSS-y}  and \ref{GSS-z}.  
              Observe that  boundary values of the coefficient $\varkappa(z)$  defined  by formula  \e{eq:wr} are given by the equalities
 \begin{equation} 
\varkappa (\lambda + i0)=
\begin{cases}
\sqrt{|\tau|} \q & \mbox{if}\q\sigma > 1, \,\tau<0, \, \lambda\in {\Bbb R}
\\ 
   \sqrt{ \lambda} \q &\mbox{if}\q\sigma < 1,   \, \lambda>0
 \\ 
  \sqrt{\lambda-\tau} \q &\mbox{if}\q\sigma = 1,    \,   \lambda>\tau 
\end{cases}
\label{eq:wr+}\end{equation}
and $\varkappa (\lambda - i0)=-\varkappa (\lambda + i0)$.

         \begin{lemma}\label{WR}
          Let one of the following three assumptions be satisfied:
        
       $1^0$ the conditions of  Theorem~\ref{GSS-x} with $\tau<0$ and $\lambda\in{\Bbb R}$ 
       
       $2^0$ the conditions of  Theorem~\ref{GSS-y} with $\gamma\lambda>0$ 
       
       $3^0$  the conditions of  Theorem~\ref{GSS-z} with $\gamma\lambda>\tau$.
       
      Then the Wronskian
       \begin{equation} 
w(\lambda) :=\frac{1}{2i}W [f(\lambda+ i0), f(\lambda- i0)]=  \gamma \varkappa (\gamma (\lambda+ i0)) > 0.
\label{eq:wr1}\end{equation}
 \end{lemma}
 

  \begin{pf} 
  Set $\varphi_{n}= \varphi_{n} (\gamma(\lambda+ i0))$, $u_{n}= u_{n} (\gamma(\lambda+ i0))$. It follows from formulas \e{eq:Gr1} and \e{eq:Jost}  that
   \[
2i w(\lambda)=-\gamma a_{n}n^{-\rho} (n+1)^{-\rho}
\Big( e^{i\varphi_{n}}e^{-i\varphi_{n+1}} u_{n} \bar{u}_{n+1}- 
e^{-i\varphi_{n}}e^{i\varphi_{n+1}} \bar{u}_{n} u_{n+1} \Big).
\]
Using condition  \e{eq:ASa}, we see that
 \begin{equation} 
w(\lambda)=- \gamma n^{\nu} \Im \big(e^{-i \theta_{n+1}}  u_{n} \bar{u}_{n+1}\big) (1+ o(1)) 
\label{eq:wr12}\end{equation}
where $\theta_{n}=\varphi_{n+1}-\varphi_{n}$ and $\nu= \sigma-2\rho$.
  Observe that
 \begin{equation} 
 \Im \big(e^{-i \theta_{n+1}}  u_{n} \bar{u}_{n+1}\big) =  \Im \big( ( u_{n} -u_{n+1})\bar{u}_{n+1}\big)- \theta_{n+1}\Re   (u_{n} \bar{u}_{n+1})+
 O(\theta_{n+1}^2).
\label{eq:wr13}\end{equation}
According to Theorems~\ref{HSS-x}, \ref{HSS-y} or \ref{HSS-z}  the first term in  the right-hand side of  \e{eq:wr13} is $O(n^{-\d +1})$ where $\d -1 > \nu$. It follows from  \e{eq:wra} that the second term  is $- \varkappa (\gamma (\lambda+ i0) ) n^{-\nu} (1+ o (1))$. Finally, the contribution of $ O(\theta_{n+1}^2)$ to \e{eq:wr12} is zero.  Therefore equality  \e{eq:wr1} is a direct consequence of \e{eq:wr12} and \e{eq:wr13}.     \end{pf}



     Let us introduce the Wronskian of the solutions $P(z)= ( P_{n} (z) )$ and $f(z)= ( f_{n} (z) )$ of equation \e{eq:Jy}:
  \begin{equation}
\Omega(z) : = W[ P (z), f (z) ] =  a_{-1} (P_{-1}(z) f_{0} (z)- P_0 (z){ f}_{-1} (z))=- a_{-1} f _{-1} (z), \q z\not\in{\cal S}.
\label{eq:J-W}\end{equation}

  \begin{lemma}\label{omega}
  The function $\Omega(z) $ is analytic in ${\Bbb C}\setminus \clos{\cal S}$ and $\Omega(z)=0$ if and only if $z$ is an eigenvalue of the operator $J$. 
In particular, $\Omega(z)\neq 0$ for $\Im z \neq 0$. 
 \end{lemma}
 
   \begin{pf}
   The analyticity of $\Omega(z) $ is a direct consequence of definition \e{eq:J-W}  because $f_{-1} (z)$, as well as all functions $f_{n} (z)$, is analytic.
   If $\Omega(z)= 0$, then $P(z)$ and $f(z)$ are proportional whence $P(z)\in \ell^2 ({\Bbb Z}_{+})$ 
by virtue of Proposition~\ref{fl2}.  Since $P_{-1} (z)=0$, it follows that $J P(z)=zP(z)$ so that
$z$ is an eigenvalue of the operator $J$. For $\Im z\neq 0$, this is impossible because $J$ is a self-adjoint  operator. Conversely, if $z$ is an eigenvalue of   $J$, then  $P(z)\in \ell^2 ({\Bbb Z}_{+})$, and hence $f(z)$ and $P(z)$  are proportional.
  \end{pf}
   
Now we are in a position to find an asymptotic behavior of the polynomials $P_{n} (\lambda)$ for $\lambda$ in the absolutely continuous spectrum (except thresholds) of the Jacobi  operator $J$.
Since the  Jost solutions $ f_{n} (\lambda\pm i0)$  are linearly independent  and $\ov{P_{n} (\lambda)} =P_{n} (\lambda)$, 
  we see that
  \begin{equation}
P_{n} (\lambda)= \ov{c(\lambda)} f_{n} (\lambda+i0)+ c(\lambda)  f_{n} (\lambda-i0)
\label{eq:J-WP}\end{equation} 
for some complex constant $c (\lambda)$. Taking the Wronskian of this equation with $ f  (\lambda+i0)$, we can  express $c(\lambda)$ via Wronskian     \e{eq:J-W}: 
  \[
- c (\lambda) W [ f (\lambda+i0),  f (\lambda-i0) ] =W[P (\lambda),  f (\lambda+ i0)]  =  \Omega (\lambda+i0) 
\]
whence 
 \[
c (\lambda)=-\frac{   \Omega (\lambda + i0)    }{ 2 i    w (\lambda)}. 
\]
In view of formula \e{eq:J-WP},  this yields the following result.

\begin{lemma}\label{HH+}
 For all $ \lambda\in \cal S $, we have the representation 
  \begin{equation}
P_{n} (\lambda)=\frac{  \Omega(\lambda- i0) f _{n} (\lambda+i0) - \Omega (\lambda + i0) f _{n} (\lambda- i0)  }{ 2 i    w (\lambda)},   \q n \in {\Bbb Z}_{+}. 
\label{eq:HH4L+}\end{equation}
 \end{lemma}
 
 Properties of  the Wronskians  $\Omega(\lambda \pm i0)$ are summarized in the following statement.
 
 \begin{theorem}\label{HX+}
 Let the assumptions of  Lemma~\ref{WR} be satisfied.
 Then the Wronskians $\Omega(\lambda + i0)$ and $\Omega (\lambda - i0)=
 \ov{ \Omega (\lambda + i0) }$  are continuous functions of $  \lambda\in \cal S $ and
 \begin{equation}
\Omega (\lambda\pm i0) \neq 0 ,\q  \lambda\in \cal S  .
\label{eq:HH5+}\end{equation}
 \end{theorem}
 
  \begin{pf}
  The functions 
 $ \Omega (\lambda\pm i0)$ are    continuous    in the same region as the Jost solutions.
  If $ \Omega (\lambda\pm i0)=0$, then according to \e{eq:HH4L+}
    $P_{n} (\lambda)=0$ for   all $n\in{\Bbb Z}_{+}$.  However, 
 $P_0 (\lambda)=1$ for all $\lambda$.
    \end{pf}

Let us set
 \begin{equation}
\kappa (\lambda) = | \Omega(\lambda+i0)|, \q
 \Omega (\lambda\pm i0) =  \kappa (\lambda) e^{\pm i \eta (\lambda) } .
\label{eq:AP+}\end{equation}
In the theory of short-range perturbations of the Schr\"odinger operator, the functions $\kappa (\lambda) $ and $\eta (\lambda)$ are known as the limit amplitude and the limit phase, respectively; the function   $\eta (\lambda)$ is also called the scattering  phase or the   phase shift.    Definition \e{eq:AP+}     fixes the phase $\eta (\lambda)$  only up to a term $2\pi m$ where $m\in{\Bbb Z}$.  We emphasize that the amplitude $\kappa (\lambda) $ and the phase $\eta (\lambda)$
 depend on the   values of the coefficients $a_{n}$ and $b_{n}$ for all $n$,  and hence they are not determined by an asymptotic behavior of $a_{n}$, $b_{n}$ as $n\to\infty$. 

Combined together,   relations \e{eq:A22G+} and \e{eq:HH4L+} yield    asymptotics of   the  orthonormal  polynomials  $P_{n} (\lambda)$.  
  

        \begin{theorem}\label{P1}
        Let one of the following three assumptions be satisfied:
        
       $1^0$ the conditions of  Theorem~\ref{GSS-x} with $\tau<0$ and $\lambda\in{\Bbb R}$ 
       
       $2^0$ the conditions of  Theorem~\ref{GSS-y} with $\gamma\lambda>0$ 
       
       $3^0$  the conditions of  Theorem~\ref{GSS-z} with $\gamma\lambda>\tau$.
       
       Let the number $\rho $ be defined by equalities \e{eq:rho}, and let $\Phi_{n}(\lambda)=\varphi_{n}(\gamma(\lambda+i0))$ where the sequences $\varphi_{n} (\lambda)$ are constructed in Theorems~\ref{GSS-x}, \ref{GSS-y}  and \ref{GSS-z}. 
       Then, for $ \lambda\in {\cal S}$,     
            \begin{equation}
P_{n}(  \lambda )  =  \kappa (\lambda) w(\lambda)^{-1} (-\gamma)^{n} n^{- \rho} \sin(  \Phi_{n} (\lambda) - \eta (\lambda)) \big(1 + o( 1)\big) , \q n\to \infty,
\label{eq:AX}\end{equation} 
where the  Wronskian $w(\lambda)$ is given by equalities   \e{eq:wr+}, \e{eq:wr1} and the amplitude $\kappa (\lambda)$ and the phase $\eta (\lambda)$ are  defined by relations \e{eq:AP+}.
 \end{theorem}
 
 We emphasize that the definitions of the numbers $\rho$ and $\Phi_{n} (\lambda)$ are different under  assumptions $1^0$, $2^0$ and $3^0$, but  relation \e{eq:AX} is true in all these cases.  Under the assumptions of Theorem~\ref{P1}  the functions  $\Phi_{n} (\lambda)$  are real and  $\Phi_{n} (\lambda)\to\infty$ so that $P_{n}(  \lambda )  $ are oscillating as $n\to\infty$.
 
    A formula completely similar to \e{eq:AX} is true for all real solutions of equation \e{eq:Jy}. Only the coefficients     $ \kappa (\lambda)$ and $ \eta (\lambda)$  are changed.
    
 \section{ Spectral results}
 
  \subsection{Resolvent. Discrete spectrum }
  


If  the minimal Jacobi operator $J_{\rm min}$  is essentially self-adjoint in the space $\ell^2 ({\Bbb Z}_{+})$, then, for $\Im z \neq 0$, equation \e{eq:Jy} has a unique (up to a constant factor) solution  $f_n(z)\in\ell^2 ({\Bbb Z}_{+})$.
    Let $I $ be the identity operator  in the space $\ell^2 ({\Bbb Z}_{+})$, and let $R(z)= (J-z I)^{-1}$ be the resolvent   of the operator
   $J=\clos J_{\rm min}$.
    Recall that the Wronskian $\Omega (z)$   of the solutions $ P_{n} (z) $ and $f_{n} (z) $ of equation \e{eq:Jy} was defined by formula \e{eq:J-W}. The following statement is very close to the corresponding result for differential operators.

 \begin{proposition}\label{res}\cite[Proposition~2.1]{inf}
 In the limit point case,
 for all $h= (h_{n})\in \ell^{2} ({\Bbb Z}_{+})$, we have
  \begin{equation}
(R (z)h)_{n} = \Omega(z)^{-1} \Big( 
f_{n} (z) \sum_{m=0}^{n} P _{m}(z)h_{m}+ P_{n} (z) \sum_{m=n+1}^{\infty}f_{m}(z) h_{m}\Big),\q \Im z \neq 0.
\label{eq:RRes}\end{equation}
 \end{proposition} 
 

 \begin{remark}\label{resA}
   Let $e_{0}, e_{1},\ldots, e_{n},\ldots$ be the canonical basis in the space $\ell^2 ({\Bbb Z}_{+})$. Then  representation  \e{eq:RRes}  can be equivalently rewritten as
\begin{equation}
\la R (z)e_{n}, e_{m}\ra =   \Omega (z)^{-1} P_{n} (z) f_{m}(z)\; \mbox{if}\; n\leq m \; \mbox{and}\; \la R (z)e_{n}, e_{m}\ra =\la R (z)e_{m}, e_{n}\ra .
\label{eq:Rrpm}\end{equation}
 \end{remark} 
 
  According to Theorem~\ref{P2+} and Proposition~\ref{D1+}, under our assumptions the   operator $J_{\rm min}$  is essentially self-adjoint in the space $\ell^2 ({\Bbb Z}_{+})$. In view of    Proposition~\ref{fl2}, in this case $f_{n}(z)$  is the Jost solution.  Thus, 
   the resolvent of the  Jacobi operator $J$  admits representation \e{eq:RRes} where  $f_{n}(z)$  is the Jost solution.



  Spectral results about  the Jacobi operators $J$  are direct consequences of representation \e{eq:RRes}. 
  As far as the discrete  spectrum is concerned, we use that
 according to Theorems~\ref{GSS-x}, \ref{GSS-y} and  \ref{GSS-z},   the functions $f_{n} (z)$, $n=-1, 0,1,\ldots$,  and, in particular,  $\Omega (z)$ are analytic     functions of $z\in{\Bbb C}\setminus \clos{\cal S}$. In view of Lemma~\ref{omega} this yields the part of Theorem~\ref{S-Adj}  concerning the discrete spectrum. Let us  state it explicitly.


   \begin{theorem}\label{Discr}
    Let  assumptions  \e{eq:ASa},  \e{eq:ASb} with $| \gamma |=1$  be satisfied. 
       
          $1^0$ If  $\sigma\in (1,3/2 ]$ and $\tau >0$,
       then the spectrum of the operator $J$ is discrete.
       
       $2^0$ 
        If  $\sigma\in (0, 1)$,
       then the spectrum of the operator $J$ is discrete on the half-axis $(-\infty, 0)$ for $\gamma=1$, and it is discrete on $( 0,\infty)$ for $\gamma= -1$.

             $3^0$   
                     If  $\sigma =1$,
       then the spectrum of the operator $J$ is discrete on the half-axis $(-\infty, \tau)$ for $\gamma=1$, and it is discrete on $( -\tau,\infty)$ for $\gamma= -1$.
 \end{theorem}
 
  \subsection{Limiting absorption principle. Continuous spectrum}
 
   Next, we consider the absolutely continuous spectrum.
   According to Theorems~\ref{GSS-x}, \ref{GSS-y} and  \ref{GSS-z},   the functions $f_{n} (z)$, $n=-1, 0,1,\ldots$,  and, in particular,  $\Omega (z)$   are continuous  up to the cut along the interval ${\cal S}$.  Therefore the following result is a direct consequence of relation \e{eq:HH5+} and representation   \e{eq:RRes}.  Recall that the set ${\cal D}\subset \ell^2({\Bbb Z}_{+})$ consists of finite lear combinations of the basis vectors $ e_{0},e_{1}, \ldots$.

 \begin{theorem}\label{AC}
   Let the assumptions of Theorems~\ref{GSS-x} for $\tau <0$, \ref{GSS-y} or \ref{GSS-z} be satisfied.  
  Then 
 for all $u,v\in {\cal D}$, the functions $\la R(z) u,v \ra$ are  continuous in $z$ up to the cut along the interval ${\cal S}$ as $z$ approaches ${\cal S}$ from upper or lower half-planes.
 \end{theorem}
 
 
 
This result   is known as the limiting absorption principle. It implies
 
 \begin{corollary}\label{Rc}
The spectrum of the operator $J$   is absolutely continuous on the closed interval $\clos{\cal S}$, except, possibly, eigenvalues at its endpoints.  In particular, it is  absolutely continuous and coincides with the whole real axis $\Bbb R$ if $\sigma\in (1,3/2]$ and $\tau<0$.
 \end{corollary}

Let us now consider the spectral projector $E(\lambda)$ of the operator $J$. By the Cauchy-Stieltjes-Privalov formula for $u,v\in {\cal D}$, its matrix elements satisfy the identity
 \begin{equation}
2\pi i \frac{d\la E (\lambda)u,v\ra} {d\lambda}=  \la R (\lambda+ i0)u,v\ra - \la R (\lambda- i0) u,v\ra ,\q \lambda\in  {\cal S}.
\label{eq:Priv}\end{equation}
Therefore, the following assertion is a direct consequence of  Theorem~\ref{AC}.

\begin{corollary}\label{RE}
 For all $u,v\in {\cal D}$, the functions $\la E(\lambda) u,v \ra $ are  continuously differentiable in $\lambda\in  {\cal S}$.
   \end{corollary}

 
   Formulas \e{eq:Rrpm} and \e{eq:Priv} allow us to calculate  the spectral family  $d E(\lambda)$ in  terms of the orthonormal polynomials and the Jost function.  Indeed, substituting the expression
\[
\la R (\lambda\pm i0)e_{n}, e_{m}\ra =  \Omega (\lambda\pm i 0)^{-1} P_{n} (\lambda) f_{m}(\lambda\pm i 0),\q n\leq m,\q \lambda\in {\cal S},
\]
 into \e{eq:Priv} and using  the identity $\Omega (\lambda-i 0) =\ov{\Omega (\lambda +i 0)} $,
 we find that
  \[
2\pi i \frac{d \la E (\lambda)e_n, e_m \ra } {d\lambda}= P_{n} (\lambda) \frac{ \Omega (\lambda-i 0)  f_{m}(\lambda+i 0)    -  \Omega (\lambda+i 0) f_{m}(\lambda-i 0)    }{|  \Omega (\lambda\pm i 0) |^2}.
\]
 Combining this representation with  formula \e{eq:HH4L+} for $P_m (\lambda) $, we obtain    the following result.

 \begin{theorem}\label{SF}
  Let the assumptions of Theorems~\ref{GSS-x} for $\tau <0$, \ref{GSS-y} or \ref{GSS-z} be satisfied.   Then   for all $n,m\in{\Bbb Z}_{+}$, we have the representation
  \begin{equation}
 \frac{d \la E (\lambda)e_n, e_m\ra } {d\lambda}= (2\pi)^{-1} w(\lambda) |\Omega (\lambda\pm i0) |^{-2} P_{n} (\lambda) P_m (\lambda)   ,\q \lambda\in \cal S,
\label{eq:sp-me}\end{equation}
where $w(\lambda)$ and $\Omega (z)$ are the Wronskians \e{eq:wr1} and \e{eq:J-W}, respectively.
 In particular,  
the spectral measure of the operator $J$ equals
  \[
d \Xi (\lambda) : =d \la E (\lambda)e_0, e_0 \ra= \xi (\lambda) d\lambda ,\q \lambda\in \cal S,
\]
where the weight  $\xi(\lambda)$ is given by the formula  
 \begin{equation}
\xi(\lambda)=  (2\pi)^{-1} w(\lambda)  |  \Omega ( \lambda\pm i0) |^{-2}   .  
\label{eq:SF1}\end{equation}
 \end{theorem} 
 
 \begin{remark}\label{stab}
  Formulas \e{eq:sp-me}, \e{eq:SF1} are also true (see \cite{inf})  in the non-critical case $|\gamma|<1$ with $w=  \sqrt{1- \gamma^2}$ and $\cal S=\Bbb R$  as well as  (see \cite{JLR})  for stabilizing coefficients satisfying \e{eq:Nev} with $w (\lambda) = 2^{-1}\sqrt{1-\lambda^2}$ and $\cal S=(-1,1)$ (if $a_{\infty}=1/2$).
     \end{remark}
 
 \begin{remark}\label{Na-Si}
  For the case $\sigma\in (0,1)$,  another representation for  the weight  $\xi(\lambda)$ was obtained in \cite{Na-Si} -- see formula (4.12) in this paper. It is difficult to compare these two representations because the Jost solutions were defined in \cite{Na-Si}  in terms of infinite products and 
  formula (4.12) contains an implicit factor (4.8).
   \end{remark}
 
 Putting together Theorem~\ref{HX+} and formula \e{eq:SF1}, we obtain
 
  \begin{theorem}\label{SFr}
  Under the assumptions of Theorem~\ref{SF}
 the weight $\xi (\lambda)$ is a continuous strictly positive function of $\lambda\in \cal S$.
 \end{theorem}

 Note that this result was deduced in  \cite{Jan-Nab-Sh} from the subordinacy theory. The assumptions of \cite{Jan-Nab-Sh} are more restrictive compared  to Theorem~\ref{SF}; in particular, it was required in  \cite{Jan-Nab-Sh} that $\sigma\in (1/2, 2/3)$.


 In view of   \e{eq:SF1} the scattering amplitude $\kappa(\lambda)$ defined by \e{eq:AP+} can be  expressed via the weight $\xi(\lambda)$:
     \[
 \kappa (\lambda) =  (2\pi)^{-1/2}   w(\lambda)^{1/2}  \xi(\lambda)^{-1/2} .
 \]
 Hence asymptotic formula   \e{eq:AX} can be rewritten as
    \[
 P_{n} (\lambda)=  \big( 2\pi w(\lambda)  \xi (\lambda)\big)^{-1/2}   (-\gamma)^{n} n^{- \rho}  \big(\sin (\Phi_{n}  (\lambda) - \eta(\lambda) ) +  o(1) \big)  
\]
as $n\to\infty$.  This form seems to be more common for the orthogonal polynomials literature.


    \end{document}